\documentclass[10pt,a4paper]{article}
\usepackage[margin=3.24cm]{geometry}
\usepackage{graphicx} 
\usepackage{amsfonts,amsmath,amssymb,bm,dsfont,mathrsfs,physics}
\usepackage{amssymb,amsthm}
\usepackage{multirow}
\usepackage{makecell}
\usepackage{enumitem,cases}
\setlist[description]{style=unboxed,leftmargin=.5em}
\setlist[itemize]{style=sameline,leftmargin=2em}
\usepackage{float}
\usepackage[dvipsnames]{xcolor}
\colorlet{xlinkcolor}{red!50!black}

\definecolor{refblue}{RGB}{26,13,171}

\usepackage{algorithm}[0.1]

\usepackage{hyperref}[6.83]
\hypersetup{
  colorlinks = true,
  allcolors = refblue,
  urlcolor = BrickRed,
}

\usepackage{cleveref}

\usepackage{subfigure}

\newtheorem{definition}{Definition}

\newtheorem{assumption}{Assumption}
\newtheorem{lemma}{Lemma}

\newtheorem{proposition}{Proposition}

\newtheorem{remark}{Remark}
\numberwithin{equation}{section}
\numberwithin{lemma}{section}
\numberwithin{example}{section}
\numberwithin{definition}{section}
\numberwithin{assumption}{section}
\numberwithin{theorem}{section}
\numberwithin{proposition}{section}
\numberwithin{corollary}{section}
\numberwithin{remark}{section}

\def\[{\begin{equation}}
\def\]{\end{equation}}
 
\setlist[description]{style=multiline,leftmargin=2em}
\setlist[itemize]{style=standard,leftmargin=2em}
\setlist[enumerate]{style=standard,leftmargin=2.2em,itemsep=3pt}

\title{\Large\bf
A Dynamical Variable-separation Method for Parameter-dependent Dynamical Systems\thanks{The research of this work was supported by the National Key R\&D Program of China (No. 2021YFA1001300), 
the National Natural Science Foundation of China (Nos. 12288201, 12271150, 12471405), 
the Youth Innovation Promotion Association of CAS, 
the Natural Science Foundation of Hunan Province (No. 2023JJ10001), 
and the Science and Technology Innovation Program of Hunan Province (No. 2022RC1190).}}

\author{Liang Chen\thanks{School of Mathematics, Hunan University, Changsha 410082, China
(\url{chl@hnu.edu.cn}, \url{yrchen@hnu.edu.cn}, \url{qli28@hnu.edu.cn}).}
\and Yaru Chen\footnotemark[2]
\and Qiuqi Li\footnotemark[2]
\and Tao Zhou\thanks{LSEC, Institute of Computational Mathematics and Scientific/Engineering Computing, Academy of Mathematics and Systems Science, Chinese Academy of Sciences, Beijing 100190, China 
 (\url{tzhou@lsec.cc.ac.cn}).}
}
\date{\today}

\begin{document}
\maketitle
 
\begin{abstract}
This paper proposes a dynamical Variable-separation method for solving parameter-dependent dynamical systems. To achieve this, we establish a dynamical low-rank approximation for the solutions of these dynamical systems by successively enriching each term in the reduced basis functions via a greedy algorithm. This enables us to reformulate the problem to two decoupled evolution equations at each enrichment step, of which one is a parameter-independent partial differential equation and the other one is a parameter-dependent ordinary differential equation. These equations are directly derived from the original dynamical system and the previously separate representation terms. Moreover, the computational process of the proposed method can be split into an offline stage, in which the reduced basis functions are constructed, and an online stage, in which the efficient low-rank representation of the solution is employed. The proposed dynamical Variable-separation method is capable of offering reduced computational complexity and enhanced efficiency compared to many existing low-rank separation techniques. Finally, we present various numerical results for linear/nonlinear parameter-dependent dynamical systems to demonstrate the effectiveness of the proposed method.
 

\medskip
\noindent 
{\bf MSCcodes:}
 35R60, 37M99, 65P99
\end{abstract}

\section{Introduction}
Dynamical systems typically model many problems in computational science and engineering, such as heat transfer, fluid dynamics, and chemically reacting flows. 
The model's input often suffers from uncertainty due to inadequate knowledge about the physical properties, geometric features, and initial or boundary conditions.
To provide accurate and reliable predictions, these uncertainties are often described by introducing a set of parameters into these systems, resulting in parameter-dependent dynamical systems.
For the numerical simulation of such problems, 
traditional approaches with high fidelity, such as the finite element method, the finite difference method, and the finite volume method, often result in millions or even billions of degrees of freedom. 
This can negatively affect the computation efficiency
and may bring challenges for numerical simulation, especially when the systems describe complex physical phenomena.

To address these problems, an efficient and practical framework is the generalized polynomial chaos (gPC) expansion \cite{Xiu2002the, Xiu2003modeling, E2008algebraic, Hoang2012REGULARITYAG, K2020long, V2020analy}, 
which approximates the solution of a parametric dynamical system in the following form
\begin{equation}
\label{gpc}
u(\bm{x},t;\bm{\xi})\approx 
u_N(\bm{x},t;\bm{\xi})
:=\sum_{i=1}^{N}\zeta_i(\bm{\xi})g_i(\bm{x},t),
\end{equation}
where $u(\bm{x},t;\bm{\xi})$ is the solution of a parametric dynamical system with 
$\bm{x}$ being the spatial variable in the physical domain $D\subset\mathbb{R}^{d}$, 
$t\in [0, T]$ being the time variable, and $\bm{\xi}$ being the parameter. 
In \eqref{gpc}, the parameter-independent spatial basis $\{g_i(\bm{x},t)\}_{i=1}^{N}$ is time-dependent, 
while the parametric coefficients $\{\zeta_i(\bm{\xi})\}_{i=1}^{N}$, often treated as elements of a Hilbert space of random functions and approximated by orthogonal polynomials, are time-independent.
These methods often require substantial prior information from the full model, which can be particularly challenging for complex dynamical systems with huge degrees of freedom. 
Another approach for approximating $u(\bm{x},t;\bm{\xi})$ is the model order reduction (MOR) \cite{T2008model,C2013nonintrusive,U2014model,P2015a,Hesthaven2022reduced,Song2024a,Li2024an}, 
which aims to approximate the solution manifold 
$$
\mathcal{M}:=\{u(\cdot,t;\bm{\xi})\mid  
t\in [0, T],\,
\bm{\xi} \in \Omega\}
$$
in a subspace with a much lower dimensionality, where $\Omega$ is a certain domain. 
The MOR methods can significantly reduce computational costs and accurately capture the most important behaviors of the original system simultaneously. 
The most prevalent MOR technique for parameter-dependent dynamical systems is the proper orthogonal decomposition (POD) method \cite{M2019data,R2019reduced,X2021non,Sleeman2022goal}, which approximates $u(\bm{x},t;\bm{\xi})$ in the form of
$$
u(\bm{x},t;\bm{\xi})\approx u_N(\bm{x},t;\bm{\xi}):=\sum_{i=1}^{N}\zeta_i(t;\bm{\xi})g_i(\bm{x}),
$$
where $\{\zeta_i(t;\bm{\xi})\}_{i=1}^{N}$ is a
basis depending on both time and parameter, 
and $\{g_i(\bm{x})\}_{i=1}^{N}$ is a spatial basis independent of time and parameter.

The success of MOR techniques is heavily based on the assumption that the solution manifold $\mathcal{M}$ has fast-decaying Kolmogorov $N$-widths \cite{A2012N-Widths,B2019Kolmo}.
This does not hold for many parameter-dependent dynamical systems, including first-order linear transport problems and hyperbolic problems \cite{M2016reduced,C2019decay,k2020model},
which makes traditional MOR techniques generally not sufficiently effective. 
One way to overcome the Kolmogorov barrier is the Dynamically Orthonormal (DO) method, as detailed in \cite{B2022breaking,D2023manifold}.
Specifically, the DO method, initiated by Sapsis and Lermusiauxthe \cite{Sapsis2009dynamical}, utilizes a more general low-rank approximation of the solution in the form of 
\begin{equation}
\label{eq-DO}
u(\bm{x},t;\bm{\xi})\approx u_N(\bm{x},t;\bm{\xi}):=\overline{u}(\bm{x},t)+\sum_{i=1}^{N}\zeta_i(t;\bm{\xi})g_i(\bm{x},t),
\end{equation}
where $\overline{u}(\bm{x},t)$ is the statistical mean field,
each $\zeta_i(t;\bm{\xi})$ is depending on both time and parameter,
and each $g_i(\bm{x},t)$ is a deterministic orthonormal field depending on time. 
This method derives an equivalent system that governs the evolution of $\overline{u}(\bm{x},t)$, $\zeta_i(t;\bm{\xi})$ and $g_i(\bm{x},t)$ directly from the original dynamical systems. 
In addition, it was studied and applied to approximate fluid equations \cite{Sapsis2012dynamical,Ueckermann2013numerical}. 
Based on the DO method, several related techniques have been developed, such as the dynamically bi-orthogonal (DyBO) method \cite{Cheng2013a1,Cheng2013a2}, the dual dynamically orthogonal formulation \cite{Musharbash2018dual}, and the dynamically/bi-orthonormal (DBO) decomposition \cite{Patil2020real-time}. 
DBO has been shown to be equivalent to the DO and DyBO methods via an invertible matrix transformation \cite{Patil2020real-time}. 
The error analysis of the DO method can be found in \cite{Musharbash2015error,Feppon2018a}, and
similar approaches with the DO method can also be found in \cite{Beck2000the,Othmar2007dynamical,Koch2007regularity}.
Furthermore, the approximation in (\ref{eq-DO}) was further advanced by Billaud-Friess and Nouy in \cite{Friess2017DynamicalMR}, in which a projection-based model order reduction method was proposed, utilizing a Galerkin projection of the full-order dynamical system onto time-dependent reduced spaces constructed by a greedy approach known as the T-greedy algorithm.

To further exploit the advantages of MOR techniques, 
in this paper, we focus on DO-type methods and develop a dynamical Variable-separation (DVS) method greedily obtaining low-rank approximations for parameter-dependent dynamical systems in the form of 
$$
u(\bm{x},t;\bm{\xi})\approx u_N(\bm{x},t;\bm{\xi}):=\sum_{i=1}^{N}\zeta_i(t;\bm{\xi})g_i(\bm{x},t),
$$
where $\{\zeta_i(t;\bm{\xi}) \}_{i=1}^{N}$ is a parameter-dependent basis and $\{ g_i(\bm{x},t)\}_{i=1}^{N}$ is a parameter-independent basis, which are all time-dependent. 
The original Variable-separation (VS) method, inspired by proper generalized decomposition \cite{Nouy2009recent,Nouy2009gener}, was proposed in \cite{Li2017a} for linear stochastic problems, and extended to stochastic saddle point problems in \cite{Jiang2018model}, to nonlinear partial differential equations (PDEs) with random inputs in \cite{Li2020a}, and to the non-overlapping domain decomposition method for stochastic PDEs in \cite{Chen2024stochastic}.
In the spirit of parametric model reduction \cite[Section 1.1]{P2015a}, the computational process of the proposed DVS method can be divided into offline and online stages,
which is the same as the VS method. 
In the offline stage, the reduced basis functions are constructed greedily. 
Specifically, at each enrichment step $k$, we solve two uncoupled evolution equations: a parameter-independent PDE and a parameter-dependent ordinary differential equation (ODE), which are derived from the original dynamical system and the previously separate representation terms. 
Firstly, the initial conditions for these subproblems are elaborately constructed. 
Subsequently, we solve the parameter-independent evolution PDE to derive the space-time basis function $g_k(\bm{x},t)$. 
Finally, the basis function $\zeta_k(t;\bm{\xi})$ depending on both time and parameter is determined by solving the parameter-dependent ODE. 
In the online stage, the approximation for a new parameter can be computed straightforwardly using the combination of these precomputed basis functions. The DVS benefits from the efficient low-rank representation of the solution, enabling rapid simulations even for a large number of parameters. 
This fulfills a critical need in design, control, optimization, and uncertainty quantification, which require repeated model evaluations over the parameter space \cite{P2015a}.
As a result, the proposed method offers reduced computational complexity and improved efficiency compared to many existing low-rank separation techniques.

The remaining parts of this paper are organized as follows. \Cref{ssec:prelim} provides the necessary notation and preliminaries.
In \Cref{sec-dvs}, we first present the general framework of the DVS method and then discuss the setting of the initial conditions for the subproblems at each enrichment step.
\Cref{sec-imple-dvs} describes the details of the DVS method for linear/nonlinear parameter-dependent dynamical systems. 
\Cref{sec-numerical examples} presents 
various numerical examples to demonstrate the performance and computational advantages of the proposed DVS method.
\Cref{sec-Conclusions} concludes this paper and offers some observations on the DVS method and its potential applications, as well as some future research topics.

\section{Problem formulation and preliminaries}
\label{ssec:prelim}
In this work, we focus on the parameter-dependent dynamical systems in the form of
\begin{equation}
\label{eq-dynamical-system}
\begin{cases}
\frac{\partial u}{\partial t}(\bm{x},t;\bm{\xi})=\mathcal{F}(u(\bm{x},t;\bm{\xi});\bm{\xi}), \ \forall ~\bm{x} \in D,~ t\in [0,T], ~ \bm{\xi} \in \Omega,\\
u(\bm{x},0;\bm{\xi})= \mu(\bm{x};\bm{\xi}),\ \forall ~\bm{x} \in D, ~ \bm{\xi} \in \Omega,\\
\mathcal{B}(u(\bm{x},t;\bm{\xi}))=g(\bm{x},t;\bm{\xi}), \ \forall ~\bm{x} \in \partial D, ~ t\in [0,T], ~ \bm{\xi} \in \Omega,
\end{cases}
\end{equation} 
where $\partial D$ is the Lipschitz continuous boundary of the domain $D$, 
$\mathcal{F}$ is a spatial differential operator, 
$\mathcal{B}$ is the boundary condition operator, $g(\cdot,t;\bm{\xi})$ is the boundary term, $\mu(\cdot;\bm{\xi})$ is the initial condition, 
and $u(\cdot,t;\bm{\xi})$ represents the solution of this parameter-dependent dynamical system, 
which is contained in a Hilbert space $\mathcal{V}$ of real-valued functions that satisfy the boundary condition with an inner product $\langle \cdot,\cdot\rangle_{\mathcal{V}}$. The associated norm is given by $\|\cdot\|_{\mathcal{V}}:=\sqrt{\langle\cdot,\cdot\rangle_{\mathcal{V}}}$.
Additionally, we define the space 
$$
\begin{array}{ll}
L^2([0,T];\mathcal{V}):=\Big\{u: [0,T]\to \mathcal{V}
\mid  \int_{0}^T\|u(t)\|_{\mathcal{V}}^2 \dd t < +\infty \Big\},
\end{array}
$$
endowed with the norm $\|u\|_{L^2([0,T];\mathcal{V})}:
=\sqrt{\int_{0}^T\|u(t)\|_{\mathcal{V}}^2 \dd t}$.
As we can see, the parameter may enter the dynamical system through the differential operator $\mathcal{F}$, the initial condition $\mu$, and/or the boundary term $g$.

For the well-posedness of (\ref{eq-dynamical-system}), we assume that the function $\mathcal{F}(u;\bm{\xi})$ is continuous with respect to $u$ and satisfies the Lipschitz condition.	
The variational formulation of (\ref{eq-dynamical-system}) reads as follows: {for all  $\bm{\xi}\in\Omega$, find $u\in L^2([0,T];\mathcal{V})$, such that}
\begin{equation}
\label{eq-varia}
\Big\langle \frac{\partial u}{\partial t}(\cdot,t;\bm{\xi}),v\Big\rangle=\langle \mathcal{F}(u(\cdot,t;\bm{\xi});\bm{\xi}),v\rangle,
\quad  \forall~ v\in \mathcal{V},~t\in [0,T],
\end{equation}
with the initial condition
\begin{equation}
\label{originalinitial}
u(\cdot ,0;\bm{\xi})=\mu(\cdot;\bm{\xi}).
\end{equation}
Here, $\langle \cdot , \cdot \rangle $ is the duality pairing between $\mathcal{V}^{*}$ (the dual space of $\mathcal{V}$) and $\mathcal{V}$.

In this paper, we aim to approximate the solution of the dynamical system (\ref{eq-dynamical-system}) with the low-rank form as
\begin{equation}
\label{eq-approx}
u(\bm{x},t;\bm{\xi})\approx u_N(\bm{x},t;\bm{\xi}):=\sum_{i=1}^{N}\zeta_i(t;\bm{\xi})g_i(\bm{x},t),
\end{equation}
where $\{\zeta_i(t;\bm{\xi}) \}_{i=1}^{N}$ is a parameter-dependent basis and $\{ g_i(\bm{x},t)\}_{i=1}^{N}$ is a parameter-independent basis, of which both are time-dependent.
Throughout this paper, we make the following blanket assumption on \eqref{eq-dynamical-system}. 
\begin{assumption}
\label{assall} 
\begin{description}
\item[(i)]
$\mathcal{F}(u;\bm{\xi})$ can be decomposed as 
$\mathcal{F}(u;\bm{\xi})=\mathcal{C}(\bm{\xi})+\mathcal{A}(u;\bm{\xi})+\mathcal{H}(u;\bm{\xi})$, 
where $\mathcal{C}:\Omega\to\mathcal{V}^*$, $\mathcal{A}:\mathcal{V}\times \Omega\to\mathcal{V}^*$ is linear about $u$, and 
$\mathcal{H}:\mathcal{V}\times \Omega\to\mathcal{V}^*$ 
is nonlinear about $u$.

\item[(ii)]
The following affine decompositions with respect to the parameter $\bm{\xi}$ hold: 
\begin{equation}
\label{eq-affine-f}
\mathcal{C}(\bm{\xi})=\sum_{i=1}^{N_C}\kappa_C^i(\bm{\xi})\mathcal{C}^i, ~
\mathcal{A}(u;\bm{\xi})=\sum_{i=1}^{N_A}\kappa_A^i(\bm{\xi})\mathcal{A}^i(u),~  
\mathcal{H}(u;\bm{\xi}) =\sum_{i=1}^{N_H}\kappa_H^i(\bm{\xi})\mathcal{H}^i(u), 
\end{equation}
where $\mathcal{C}^i\in\mathcal{V}^*$, $\mathcal{A}^i:\mathcal{V}\to\mathcal{V}^*$, $\mathcal{H}^i:\mathcal{V}\to\mathcal{V}^*$, and the scalar functions $\kappa_C^i(\bm{\xi})$, $\kappa_A^i(\bm{\xi})$, $\kappa_H^i(\bm{\xi}):\Omega\to \mathbb{R}$.

\item[(iii)]
The initial condition of \eqref{eq-dynamical-system} is given by
\begin{equation}
\label{eq-u0-affine}
{\mu(\bm{x};\bm{\xi})}=\sum_{i=1}^{N_{t_0}}p^i(\bm{\xi})q^i(\bm{x}),\ \forall ~\bm{\xi}\in\Omega,~\bm{x} \in D,
\end{equation}
where $\{p^i(\bm{\xi})\}_{i=1}^{N_{t_0}}$ are $\bm{\xi}$-dependent functions, $\{q^i(\bm{x})\}_{i=1}^{N_{t_0}}$ are independent of the parameter $\bm{\xi}$. 
\end{description}
\end{assumption}

Here we should emphasize that the affine decompositions are crucial to achieving the online-offline decomposition, 
and \Cref{assall} is very moderate and reasonable.  
Specifically, if it is not satisfied, we can use the VS method for multivariate functions in \cite{Li2017a} to obtain such an affine expansion approximation with negligible inaccuracy.

\section{A dynamical Variable-separation method}
\label{sec-dvs}
This section develops the DVS method for \eqref{eq-dynamical-system}. 
We first provide a general framework for this method,
where the basis functions are generated by a greedy algorithm.
Then we discuss the setting of initial conditions for subproblems at each enrichment step, which is one of the most essential ingredients for the proposed method.

\subsection{A general framework of the DVS method}
\label{sec-dvs-framework}
We start by presenting the general framework from the offline stage of the DVS method,
which constructs the basis functions $\{g_i(\bm{x},t)\}_{i=1}^{N}$ and $\{\zeta_i(t;\bm{\xi})\}_{i=1}^{N}$ in the approximate solution in (\ref{eq-approx}) of the dynamical system (\ref{eq-dynamical-system}).

Let $\Xi\subset \Omega$ be the training set consisting of a finite number of samples, and
${u}(\bm{x},t;\bm{\xi})$ be the unknown solution of \eqref{eq-dynamical-system} we want to approximate. 
Firstly, we set the step counter $k$ to $1$ and arbitrarily choose $\bm{\xi}_1$ in $\Xi$.
Then update $\Xi:=\Xi \textbf{\textbackslash} \bm{\xi}_1$. 
When $g_1(\bm{x},t)$ is taken as the solution of (\ref{eq-varia}) under the initial condition \eqref{originalinitial} with $\bm{\xi}=\bm{\xi}_1$, 
one can let $\zeta_1(t;\bm{\xi})$ be an unknown basis function such that 
$u_1(\bm{x},t;\bm{\xi})=g_1(\bm{x},t)\zeta_1(t;\bm{\xi})$. 
Taking $u=u_1$ and $v=g_1(\bm{x},t)$ in \eqref{eq-varia} one gets
\begin{equation}
\label{eq-frame-1}
\Big\langle \frac{\partial (g_1(\bm{x},t)\zeta_1(t;\bm{\xi}))}{\partial t},g_1(\bm{x},t)\Big\rangle=\langle \mathcal{F}(g_1(\bm{x},t)\zeta_1(t;\bm{\xi});\bm{\xi}),g_1(\bm{x},t)\rangle.
\end{equation}
Then by solving the above parameter-dependent ODE under a certain initial condition (which will be given by \eqref{initialconditionfirst}), one gets the basis function $\zeta_1(t;\bm{\xi})$.

At the $k$-th step ($k\ge 2$) of the offline stage, by denoting ${e}(\bm{x},t;\bm{\xi}):={u}(\bm{x},t;\bm{\xi})-{u}_{k-1}(\bm{x},t;\bm{\xi})$, one can set
$$
\bm{\xi}_k\in \mathop{\text{argmax}}\limits_{\bm{\xi}\in\Xi} \bigtriangleup_k(\bm{\xi}),
$$
where the value of $\bigtriangleup_k(\bm{\xi})$ is a proper upper bound estimation of 
$||e(\bm{x},t;\bm{\xi})||_{L^2([0, T];\mathcal{V})}$ (c.f.  \Cref{sec-error} for the details), or merely itself if the computation is affordable. 
When the error estimator $\bigtriangleup_k(\bm{\xi}_k)$ is sufficiently small, the offline stage is stopped, or else we update $\Xi:=\Xi \textbf{\textbackslash} \bm{\xi}_k$ and perform the following procedure. 
For convenience, we denote $e=u-u_{k-1}$ and write \eqref{eq-varia} as 
\begin{equation}
\label{eq-frame-pde}
~~~\Big\langle \frac{\partial (e+u_{k-1})}{\partial t}(\bm{x},t;\bm{\xi}),v\Big\rangle
=\langle \mathcal{F}((e+u_{k-1})(\bm{x},t;\bm{\xi});\bm{\xi}),v\rangle,
\quad \forall~ v\in \mathcal{V},~t\in [0,T].
\end{equation}
Taking $g_k(\bm{x},t)$ as the solution of (\ref{eq-frame-pde}) at $\bm{\xi}=\bm{\xi}_k$ under a certain initial condition that will be detailed later,   
one can let $\zeta_k(t;\bm{\xi})$ be a unknown basis function and define 
${\tilde e}(\bm{x},t;\bm{\xi})=g_k(\bm{x},t)\zeta_k(t;\bm{\xi})$,
Then, under a certain initial condition (which will be specified in the next subsection in detail), one can determine $\zeta_k(t;\bm{\xi})$ by taking $e={\tilde e}$ and $v=g_k(\bm{x},t)$ in (\ref{eq-frame-pde}), i.e., 
solving the equation
\begin{equation}
\begin{aligned}
\label{eq-frame-rode}
&\Big\langle \frac{\partial (g_k(\bm{x},t)\zeta_k(t;\bm{\xi})+u_{k-1}(\bm{x},t;\bm{\xi}))}{\partial t},g_k(\bm{x},t)\Big\rangle \\
&=\langle \mathcal{F}(g_k(\bm{x},t)\zeta_k(t;\bm{\xi})+u_{k-1}(\bm{x},t;\bm{\xi});\bm{\xi}),g_k(\bm{x},t)\rangle.
\end{aligned}
\end{equation}

Note that the two subproblems at each enrich step (of the offline stage) need proper initial conditions to ensure the accuracy of the approximation. This will be discussed in detail in the next subsection.

To show the DVS method for parameter-dependent dynamical systems more clearly, we will present two examples that fit this abstract framework, that is, (i) linear parameter-dependent dynamical systems; (ii) nonlinear parameter-dependent dynamical systems (taking Burgers equation as an example) in \Cref{sec-imple-dvs}.

\subsection{Initial conditions for subproblems}
To facilitate comprehension and implementation, we will discuss the details of the initial conditions for the subproblems, that is, the parameter-independent PDEs \eqref{eq-frame-pde} and the parameter-dependent ODEs \eqref{eq-frame-1} and \eqref{eq-frame-rode}, during the offline stage.

Recall that during the offline stage of the DVS method, we took the solution of (\ref{eq-varia}) with $\bm{\xi}=\bm{\xi}_1$ as $g_1(\bm{x},t)$ in the first step.
Therefore, we have $g_1(\bm{x},0)=\mu(\bm{x};\bm{\xi}_1)$.
To get the initial condition $\zeta_1(0;\bm{\xi})$ of the parameter-dependent ODE (\ref{eq-frame-1}), we approximate the initial condition \eqref{originalinitial} via $u_1(\cdot,0;\bm{\xi})\approx \mu (\cdot;\bm{\xi})$
with $u_1$ taking the following form
$$
u_1(\bm{x},0; \bm{\xi})
=\mu(\bm{x};\bm{\xi}_1)\zeta_1(0;\bm{\xi})=g_1(\bm{x},0)\zeta_1(0;\bm{\xi}),
$$
in the sense that 
$\langle u_1(\bm{x},0;\bm{\xi}), g_{1}(\bm{x},0)\rangle
=\langle \mu(\bm{x}; \bm{\xi}), g_{1}(\bm{x},0)\rangle$ holds. 
In the following discussion, we simplify the notation of the functions at $t=0$, e.g., we denote $g_1(\bm{x},0)$ as $g_{1,0}(\bm{x})$. 
Then, by using the separable formulation of the initial condition (\ref{eq-u0-affine}) in \Cref{assall}, 
we obtain that  
$$
\Big\langle \sum_{i=1}^{N_{t_0}}p^i(\bm{\xi})q^i(\bm{x}),g_{1,0}(\bm{x}) \Big\rangle = \big \langle g_{1,0}(\bm{x}),g_{1,0}(\bm{x}) \big \rangle\zeta_{1,0}(\bm{\xi}).
$$
This implies 
\begin{equation}
\label{initialconditionfirst}
\zeta_{1,0}(\bm{\xi})=\sum_{i=1}^{N_{t_0}}  \frac{\langle q^i(\bm{x}),g_{1,0}(\bm{x}) \rangle}{\langle g_{1,0}(\bm{x}),g_{1,0}(\bm{x}) \rangle}p^i(\bm{\xi}),
\end{equation}
which constitutes an initial condition for \eqref{eq-frame-1}.

For the subsequent iteration steps with $k\ge2$, one can see that 
the error $e(\bm{x},t;\bm{\xi})$ at $t=0$ is given by
\begin{equation}
\label{eq-ingk}
e_0(\bm{x};\bm{\xi}) =\mu(\bm{x};\bm{\xi})-u_{k-1,0}(\bm{x};\bm{\xi})
=\mu(\bm{x};\bm{\xi})-\sum_{j=1}^{k-1}g_{j,0}(\bm{x})\zeta_{j,0}(\bm{\xi}).
\end{equation}
Then we get the initial condition of the parameter-independent residual equation (\ref{eq-frame-pde}) with $\bm{\xi}=\bm{\xi}_k$ immediately from (\ref{eq-ingk}), i.e.,
\begin{equation}
\label{eq-initial-deter-residual}
g_{k,0}(\bm{x})=e_{0}(\bm{x};\bm{\xi}_k)=\mu(\bm{x};\bm{\xi}_k)-\sum_{j=1}^{k-1}g_{j,0}(\bm{x})\zeta_{j,0}(\bm{\xi}_k).
\end{equation}
Then, using $\tilde e_{0}(\bm{x};\bm{\xi}):=g_{k,0}(\bm{x})\zeta_{k,0}(\bm{\xi})$ as an approximation of $e_0(\bm{x};\bm{\xi})$ in the sense that
$$\langle \tilde e_0(\bm{x};\bm{\xi}), g_{k,0}(\bm{x})\rangle=\langle e_0(\bm{x};\bm{\xi}), g_{k,0}(\bm{x})\rangle, $$ 
one can call \eqref{eq-ingk} to obtain
\begin{equation}
\label{eq-intial-rode}
\begin{aligned}
\zeta_{k,0}(\bm{\xi})
&
=\frac{\langle e_{0}(\bm{x};\bm{\xi}), g_{k,0}(\bm{x}) \rangle}{\langle g_{k,0}(\bm{x}), g_{k,0}(\bm{x}) \rangle}
=\frac{\langle \mu(\bm{x};\bm{\xi})-\sum\limits_{j=1}^{k-1}g_{j,0}(\bm{x})\zeta_{j,0}(\bm{\xi}),g_{k,0}(\bm{x})\rangle}{\langle g_{k,0}(\bm{x}), g_{k,0}(\bm{x}) \rangle}
\\[2mm]
& 
=\sum\limits_{i=1}^{N_{t_0}}\frac{\langle q^i(\bm{x}),g_{k,0}(\bm{x})\rangle}{\langle g_{k,0}(\bm{x}), g_{k,0}(\bm{x}) \rangle}p^i(\bm{\xi})
-\sum\limits_{j=1}^{k-1}\frac{\langle g_{j,0}(\bm{x}),g_{k,0}(\bm{x})\rangle}{\langle g_{k,0}(\bm{x}), g_{k,0}(\bm{x}) \rangle}\zeta_{j,0}(\bm{\xi}),
\end{aligned}
\end{equation}
which constitutes an initial condition for \eqref{eq-frame-rode}.

\begin{remark}
Note that if the initial condition of the original dynamical system satisfies $\mu(\bm{x};\bm{\xi})=0$, we can naturally define the initial conditions for the subproblems as $g_{k,0}(\bm{x})=0$ and $\zeta_{k,0}(\bm{\xi})=0$ in each enrichment step $k$.
\end{remark}

Now, based on the above initial conditions, we can present the offline stage in the DVS method as \Cref{algorithm-dvs}.
\begin{algorithm}[H]
\caption{The offline stage of the DVS method\label{algorithm-dvs}}
\textbf{Input:} {The system (\ref{eq-dynamical-system}), the training set $\Xi\subset \Omega$, and the error tolerance $\varepsilon>0$.}

\textbf{Output:} {Surrogate model $u_N(\bm{x},t;\bm{\xi})= \sum_{i=1}^{N}\zeta_i(t;\bm{\xi})g_i(\bm{x},t)$.} 

1: Initialize the iteration counter $k=1$, arbitrarily choose $\bm{\xi}_1\in \Xi$;

$~~$1.1: Calculate $g_1(\bm{x},t)$ by solving (\ref{eq-varia})
under the initial condition \eqref{originalinitial} with $\bm{\xi}=\bm{\xi}_1$;

$~~$1.2: Calculate $\zeta_1(t;\bm{\xi})$ by solving (\ref{eq-frame-1}) under the initial condition  \eqref{initialconditionfirst};

$~~$1.3: Update $\Xi:=\Xi \textbf{\textbackslash} \bm{\xi}_1$ and set $k=2$;

2: If $\Xi=\emptyset$ or $\max\{\bigtriangleup_k(\bm{\xi})\mid {\bm{\xi}\in\Xi} \}<\varepsilon$, set $N=k-1$ and \textbf{terminate}, else \textbf{continue};

3: choose $\bm{\xi}_k \in\mathop{\text{argmax}}\{\bigtriangleup_k(\bm{\xi})\mid \bm{\xi}\in\Xi\}$;

$~~$3.1: Calculate $g_k(\bm{x},t)$ by solving (\ref{eq-frame-pde}) with $\bm{\xi}=\bm{\xi}_k$ under the initial condition \eqref{eq-initial-deter-residual};

$~~$3.2: Calculate $\zeta_k(t;\bm{\xi})$ by solving (\ref{eq-frame-rode}) under the initial condition   \eqref{eq-intial-rode}; 

$~~$3.3: Update $\Xi:=\Xi \textbf{\textbackslash} \bm{\xi}_k$, set $k:=k+1$, and \textbf{return} to Step 2.
\end{algorithm}
The online stage of the DVS method for parameter-dependent dynamical systems is described in \Cref{algorithm-dvs-online}. 
\begin{algorithm}[H]
\caption{The online stage of the DVS method \label{algorithm-dvs-online}}
\textbf{Input:} {A parameter $\bm{\bar\xi}\in\Omega$.}

\textbf{Output:} {The approximated solution ${\hat u}(\cdot, \cdot ; \bm{\bar\xi})$.}

1. Use the output of \Cref{algorithm-dvs} to calculate ${\hat u}(\cdot,\cdot ;\bm{\bar\xi}):=u_N(\cdot,\cdot ; \bm{\bar\xi} )$.
\end{algorithm}

Algorithms \ref{algorithm-dvs} and \ref{algorithm-dvs-online} constitute the whole framework of the DVS method that we propose in this paper. 
Whenever the surrogate model is obtained from the former, one can calculate the online solutions via an explicit expression, which is extraordinarily simple.
However, in practical implementations, this is realistic only for a very limited class of simple problems for which the corresponding parameter-dependent ODEs (equations \eqref{eq-frame-1} and \eqref{eq-frame-rode}) can be analytically solved. 
Specifically, for the vast majority of problems, explicit analytical expressions of the parameter-dependent basis \(\zeta_k(t;\bm{\xi})\), which should be obtained from Step 1.2 and Step 3.2 of \Cref{algorithm-dvs}, are not available, and only numerical solutions can be explicitly expressed in discrete time instances.
As can be observed in the next section, the procedure of generating the \(\zeta_k(t;\bm{\xi})\) includes complicated combinations of functions of $\bm{\xi}$.
In practice, it is often more favorable to store a few scalars involved in the expression of \(\zeta_k(t;\bm{\xi})\) instead of the final form of this basis and recompute it in the online stage for the specific instance $\bm{\bar\xi}$ only by solving very simple ODEs.
In particular, this recomputation is easier to implement and highly efficient, incurs negligible computational overhead, and operates independently of the full-order model; details will be described in \Cref{sec-imple-dvs}.

\section{Implementation of the DVS method}
\label{sec-imple-dvs}
In this section, to illustrate the DVS method more clearly, we present its implementations on two classes of concrete problems.
One is the linear parameter-dependent dynamical system and the other is the nonlinear parameter-dependent dynamical system represented by the Burgers equation.

\subsection{The DVS method for linear dynamical systems}
\label{sec-dvs-linear}
Firstly, under \Cref{assall}, we consider the more specific linear case that
$$
\mathcal{F}(u;\bm{\xi})=\mathcal{C}(\bm{\xi})+\mathcal{A}(u;\bm{\xi}).
$$
In this case, \eqref{eq-varia} in the variational formulation becomes 
$$
\Big\langle \frac{\partial u}{\partial t}(\cdot,t;\bm{\xi}),v\Big\rangle=\langle \mathcal{C}(\bm{\xi})+\mathcal{A}(u(\cdot,t;\bm{\xi});\bm{\xi}),v\rangle,
\quad 
\forall~ v\in \mathcal{V},~ t\in [0,T].
$$
As a result, at the $k$-th iteration ($k\ge 2$) of  \Cref{algorithm-dvs},  equation \eqref{eq-frame-pde} turns to 
\begin{equation}
    \label{eq-linear-residual}
    \Big\langle \frac{\partial e}{\partial t}(\bm{x},t;\bm{\xi}),v\Big\rangle
    -\langle \mathcal{A}(e(\bm{x},t;\bm{\xi});\bm{\xi}),v\rangle=\langle r_k(\bm{x},t;\bm{\xi}),v\rangle,
\end{equation}
where the residual is defined by 
$$
r_k(\bm{x},t;\bm{\xi}):=
\mathcal{A}(u_{k-1}(\bm{x},t;\bm{\xi});\bm{\xi})- \frac{\partial u_{k-1}}{\partial t}(\bm{x},t;\bm{\xi})+\mathcal{C}(\bm{\xi}).
$$

After calculating $g_k(\bm{x},t)$ by solving (\ref{eq-linear-residual}) with $\bm{\xi}=\bm{\xi}_k$ under the initial condition \eqref{eq-initial-deter-residual}, one concentrates on obtaining the parameter-dependent basis function 
 $\zeta_k(t;\bm{\xi})$ from 
 \eqref{eq-frame-rode}, which can be written as  
\begin{equation}
\label{halfreformulation}
\begin{aligned}
\Big\langle \frac{\partial (g_k(\bm{x},t)\zeta_k(t;\bm{\xi}))}{\partial t}, g_k(\bm{x},t)\Big\rangle -  
\langle \mathcal{A}(g_k(\bm{x},t)\zeta_k(t;\bm{\xi});\bm{\xi}),g_k(\bm{x},t)\rangle\\
=\langle r_k(\bm{x},t;\bm{\xi}),&g_k(\bm{x},t)\rangle.
\end{aligned}
\end{equation}
In fact, for the first iteration of \Cref{algorithm-dvs} with
$k=1$, by taking $e(\bm{x},t;\bm{\xi})=u(\bm{x},t;\bm{\xi})$ (i.e., define $e:=u-u_{k-1}$ with $u_{k-1}\equiv 0$ for $k=1$) and $r_1(\bm{x},t;\bm{\xi}) \equiv \mathcal{C}(\bm{\xi})$, the above equation is also a reformulation of \eqref{eq-frame-1}.

One can omit the variable $\bm{x}$ in the inner products above and use $(g_k)_t(t)$ to denote the partial differential of $g_k(t)$ with respect to the time-variable $t$ for convenience, i.e., $\frac{\partial g_k}{\partial t}(t)$ (such a strategy will be used to simplify notation in the following discussion).
Then, the equation  \eqref{halfreformulation} can be equivalently reformulated as
\begin{equation}
    \label{eq-linear-rode}
    \begin{aligned}
        \langle (g_k)_t(t),g_k(t)\rangle \zeta_k(t;\bm{\xi})&+
        \langle g_k(t),g_k(t)\rangle (\zeta_k)_t(t;\bm{\xi})
        \\&-\langle \mathcal{A}(g_k(t);\bm{\xi}),g_k(t)\rangle \zeta_k(t;\bm{\xi})=\langle r_k(t;\bm{\xi}),g_k(t)\rangle.
    \end{aligned}
\end{equation}

Now we discuss the details to construct the parameter-dependent basis function by solving the parameter-dependent ODE (\ref{eq-linear-rode}).  Our goal is to obtain the explicit expression of the solution $\zeta_k(t;\bm{\xi})$ in a fixed time $t$, which is a key component to achieving efficient online computing.
By invoking the affine expansion in \eqref{eq-affine-f} for $\mathcal{A}(u;\bm{\xi})$ and $\mathcal{C}(\bm{\xi})$ and the separate formulation of $u_{k-1}$ obtained from the previous iteration, the equation (\ref{eq-linear-rode}) can be rewritten as 
\begin{equation}
\label{eq-e-affine}
\begin{aligned}
&\langle (g_k)_t(t),g_k(t)\rangle \zeta_k(t;\bm{\xi})
\\
&\qquad+ 
\langle g_k(t),g_k(t)\rangle (\zeta_k)_t(t;\bm{\xi})-\zeta_k(t;\bm{\xi})\sum_{i=1}^{N_A}\kappa_A^i(\bm{\xi}) \langle \mathcal{A}^i(g_k(t)),g_k(t)\rangle\\
=&\Big\langle\sum_{i=1}^{N_A}\sum_{j=1}^{k-1}\kappa_A^i(\bm{\xi})\zeta_{j}(t;\bm{\xi}) \mathcal{A}^i(g_{j}(t))- \sum_{j=1}^{k-1}(\zeta_{j}(t;\bm{\xi})g_j(t))_t,g_k(t)\Big\rangle
+\sum_{i=1}^{N_C}\kappa_C^i(\bm{\xi})\langle \mathcal{C}^i,g_k(t)\rangle.
\end{aligned}
\end{equation}
In practical implementation, we divide the time interval $[0,T]$ into $N_t$ sub-intervals with equal step size $\tau=T/N_t$, and use $g_{k,n}$ and $\zeta_{k,n}(\bm{\xi})$ to denote $g_k(t_n)$ and $\zeta_k(t_n;\bm{\xi})$ respectively. 
Furthermore, we use first-order differences to approximate first-order derivatives as follows
$$
(g_k)_t(t_{n+1})\approx\frac{g_{k,n+1}-g_{k,n}}{\tau},
\quad 
(\zeta_k)_t(t_{n+1};\bm{\xi})\approx\frac{\zeta_{k,n+1}(\bm{\xi})-\zeta_{k,n}(\bm{\xi})}{\tau}.
$$
Invoking the backward Euler method, \eqref{eq-e-affine} can be discretized as 
$$
    \begin{aligned}
        &\Big\langle \frac{g_{k,n+1}-g_{k,n}}{\tau},g_{k,n+1}\Big\rangle \zeta_{k,n+1}(\bm{\xi})+
        \langle g_{k,n+1},g_{k,n+1}\rangle \frac{\zeta_{k,n+1}(\bm{\xi})-\zeta_{k,n}(\bm{\xi})}{\tau}\\
        -&\zeta_{k,n+1}(\bm{\xi})\sum_{i=1}^{N_A}\kappa_A^i(\bm{\xi}) \langle \mathcal{A}^i(g_{k,n+1}),g_{k,n+1}\rangle\\
        =&\sum_{i=1}^{N_A} \sum_{j=1}^{k-1}\kappa_A^i(\bm{\xi})\zeta_{j,n+1}(\bm{\xi})\langle \mathcal{A}^i(g_{j,n+1}),g_{k,n+1}\rangle
        +\sum_{i=1}^{N_C}\kappa_C^i(\bm{\xi})\langle \mathcal{C}^i,g_{k,n+1}\rangle
        \\
        -&\sum_{j=1}^{k-1}\frac{\zeta_{j,n+1}(\bm{\xi})-\zeta_{j,n}(\bm{\xi})}{\tau}\langle g_{j,n+1},g_{k,n+1}\rangle
        -\sum_{j=1}^{k-1}\zeta_{j,n+1}(\bm{\xi})\Big\langle 
        \frac{g_{j,n+1}-g_{j,n}}{\tau},g_{k,n+1}\Big\rangle.
    \end{aligned}
$$
Then we can get the explicit expression of $\zeta_{k,n+1}(\bm{\xi})$ as 
\begin{equation}
    \label{eq-iter-t}
    \begin{aligned}
        \zeta_{k,n+1}(\bm{\xi})=\frac{c_{n+1}\zeta_{k,n}(\bm{\xi})+s_{n+1}(\bm{\xi})}{l_{n+1}(\bm{\xi})},
        \quad 
        n=0,\ldots, N_t-1,
    \end{aligned}
\end{equation}
where
$$
\begin{array}{lll}
\displaystyle 
c_{n+1}=\frac{\langle g_{k,n+1},g_{k,n+1}\rangle}{\tau},
\\[2mm]
\displaystyle
l_{n+1}(\bm{\xi}) =2c_{n+1}-\frac{\langle g_{k,n},g_{k,n+1}\rangle}{\tau}
        -\sum_{i=1}^{N_A}\kappa_A^i(\bm{\xi}) \langle \mathcal{A}^i(g_{k,n+1}),g_{k,n+1}\rangle,\\
\displaystyle
s_{n+1}(\bm{\xi})
=\sum_{i=1}^{N_A} \sum_{j=1}^{k-1}\kappa_A^i(\bm{\xi})\zeta_{j,n+1}(\bm{\xi})\langle \mathcal{A}^i(g_{j,n+1}),g_{k,n+1}\rangle
+\sum_{i=1}^{N_C}\kappa_C^i(\bm{\xi})\langle \mathcal{C}^i,g_{k,n+1}\rangle
\\
\quad  
\displaystyle
-\sum_{j=1}^{k-1}\frac{\zeta_{j,n+1}(\bm{\xi})-\zeta_{j,n}(\bm{\xi})}{\tau}\langle g_{j,n+1},g_{k,n+1}\rangle
-\sum_{j=1}^{k-1}\zeta_{j,n+1}(\bm{\xi})\Big\langle 
\frac{g_{j,n+1}-g_{j,n}}{\tau},g_{k,n+1}\Big\rangle.
\end{array}
$$

Now, the numerical solution of the parameter-dependent ODE (\ref{eq-linear-rode}) can be explicitly given by (\ref{eq-iter-t}). 
Recall that (\ref{eq-linear-rode}) is exactly the equation (\ref{eq-frame-1}) for $k=1$ or \eqref{eq-frame-rode} for $k\geq 2$ in the general framework of the DVS method.
Therefore, thanks to (\ref{eq-iter-t}), one can obtain the explicit affine expression of $\zeta_{k,n+1}(\bm{\xi})$, $n=0,\ldots, N_t-1$.
Generally, a sufficiently large $N_t$ should be specified to ensure the accuracy of the discretized model involving \eqref{eq-iter-t}. This presents challenges to both storage and computation.  
Fortunately, in practice, we can address this issue by storing only a few scalars, i.e.,
\begin{equation}
\label{data1}
\begin{cases}
\langle g_{k,n+1},g_{k,n+1}\rangle, 
\ 
\langle g_{k,n},g_{k,n+1}\rangle, 
\\
\langle\mathcal{A}^i(g_{k,n+1}),g_{k,n+1}\rangle,
\langle \mathcal{A}^i(g_{j,n+1}),g_{k,n+1}\rangle,
&i=1,\dots, N_A,
\\
\langle g_{j,n+1},g_{k,n+1}\rangle, 
\langle g_{j,n},g_{k,n+1}\rangle, 
&
j=1, \dots,k-1,
\\
\langle \mathcal{C}^i,g_{k,n+1}\rangle,
& i=1, \dots, N_C, 
\end{cases}
\end{equation}
instead of the data of $\zeta_{k,n+1}(\bm{\xi})$ ($n=0,\ldots, N_t-1$). 
Thus, we obtained from the offline stage not a complete surrogate model of the dynamical system but an intermediate form. 
Specifically, in the online stage, to obtain the approximate solution $u(\bm{x},t;\bm{\bar\xi})$ for a new parameter $\bm{\bar\xi}$, the spatial basis $\{g_k(\bm{x},t)\}_{k=1}^{N}$ obtained in the offline stage can be used directly, 
but we need to calculate the value of $\zeta_k(t;\bm{\bar\xi}), k=1,\cdots, N$.
To achieve this, we can use \eqref{eq-e-affine} by fixing $\bm{\xi}\equiv\bm{\bar\xi}$, i.e., the equation obtained the parametric basis in the offline stage, which can be rewritten as 
\begin{equation}
\label{eq-online}
\alpha(t)(\zeta_k)_t(t;\bm{\bar\xi}) + \beta(t;\bm{\bar\xi}) \zeta_k(t;\bm{\bar\xi}) = \gamma(t;\bm{\bar\xi}),  
\end{equation}
where 
$$
\begin{cases}
\alpha(t) = \langle g_k(t),g_k(t)\rangle,
\\[1mm]
\beta(t;\bm{\bar\xi})=\langle (g_k)_t(t),g_k(t)\rangle -\sum_{i=1}^{N_A}\kappa_A^i(\bm{\bar\xi}) \langle \mathcal{A}^i(g_k(t)),g_k(t)\rangle,
\\[1mm]
\gamma(t;\bm{\bar\xi}) = \sum_{i=1}^{N_A}\sum_{j=1}^{k-1}\kappa_A^i(\bm{\bar\xi}) \zeta_j(t;\bm{\bar\xi})\langle\mathcal{A}^i(g_j(t)),g_k(t)\rangle 
+ \sum_{i=1}^{N_C}\kappa_C^i(\bm{\bar\xi})\langle \mathcal{C}^i,g_k(t)\rangle
\\[1mm]
\qquad 
\qquad - \sum_{j=1}^{k-1} (\zeta_j(t;\bm{\bar\xi}))_t \langle g_j(t),g_k(t)\rangle
- \sum_{j=1}^{k-1} \zeta_j(t;\bm{\bar\xi}) \langle (g_j(t))_t,g_k(t)\rangle,
\end{cases}
$$
in which the terms independent of the parameter were given by \eqref{data1}, which were stored in the offline stage. 
In the online stage, we can obtain the value of the parametric basis functions by solving the simple ODE \eqref{eq-online} since it is independent of the spatial discretization. 
In summary, 
the online stage is independent of the full-order mode and uses a computationally more economical approach to avoid using a large amount of memory, but incurs negligible computational overhead.

Finally, we should mention that both $\beta(t;\bm{\xi})$ and $\gamma(t;\bm{\xi})$ are affine with respect to the parameter. Therefore, we can achieve efficient computation in the online stage even for a large number of parameter samples.

\subsection{The DVS method for nonlinear dynamical systems: illustration via Burgers equation}
\label{sec-dvs-burgers}
In this subsection, we present the implementation of the DVS method for nonlinear dynamical systems. We take the parameter-dependent Burgers equation as an illustrative example and mention that a similar procedure can be derived for other equations (such as the Allen-Cahn equation in \cref{num-ac}).
In this case, one has
$$
\mathcal{F}(u;\bm{\xi})=\mathcal{A}(u;\bm{\xi})+\mathcal{H}(u;\bm{\xi}),
$$
where $\mathcal{A}(u;\bm{\xi})=\kappa(\bm{\xi})\frac{\partial^2 u}{\partial \bm{x}^2}$ and $\mathcal{H}(u;\bm{\xi})=-u\frac{\partial u}{\partial \bm{x}}$. 
So the variational formulation \eqref{eq-varia} for the parameter-dependent Burgers equation becomes
$$
\Big\langle \frac{\partial u}{\partial t}(\cdot,t;\bm{\xi}),v\Big\rangle =\langle \mathcal{A}(u(\cdot,t;\bm{\xi});\bm{\xi}) + \mathcal{H}(u(\cdot,t;\bm{\xi});\bm{\xi}) ,v\rangle,
\quad \forall~ v\in \mathcal{V},~ t\in [0,T].
$$
In the $k$-th iteration ($k\ge 2$) of  \Cref{algorithm-dvs}, one can define the function 
$$
\mathcal{L}(e,u_{k-1}): =e\frac{\partial u_{k-1}}{\partial \bm{x}}+ u_{k-1}\frac{\partial e}{\partial \bm{x}}- \mathcal{A}(e;\bm{\xi}),
$$
which is linear in $e$ and $u_{k-1}$. Then the equation in \eqref{eq-frame-pde} turns to
\begin{equation}
\label{eq-burgers-pde}
\begin{aligned}
&\Big\langle \frac{\partial e}{\partial t}(\bm{x},t;\bm{\xi}),v\Big\rangle
+\Big\langle e\frac{\partial e}{\partial \bm{x}}(\bm{x},t;\bm{\xi}),v\Big\rangle\\
&\qquad\qquad +\left\langle \mathcal{L}(e(\bm{x},t;\bm{\xi}),u_{k-1}(\bm{x},t;\bm{\xi})),v\right\rangle
=\langle r_k(\bm{x},t;\bm{\xi}),v\rangle,
\end{aligned}
\end{equation}
where the residual is defined by
$$
r_k(\bm{x},t;\bm{\xi}):=\mathcal{A} (u_{k-1}(\bm{x},t;\bm{\xi});\bm{\xi}) - u_{k-1}\frac{\partial u_{k-1}}{\partial \bm{x}}(\bm{x},t;\bm{\xi})-\frac{\partial u_{k-1}}{\partial t}(\bm{x},t;\bm{\xi}).
$$

By solving (\ref{eq-burgers-pde}) with fixing $\bm{\xi} =\bm{\xi}_k$ under the initial condition (\ref{eq-initial-deter-residual}), we can obtain $g_k(\bm{x},t)$.
In practical simulation, applying the backward Euler method for time discretization to \eqref{eq-burgers-pde} obtains
$$
\begin{aligned}
&\Big\langle 
\frac{e_{n+1}(\bm{x};\bm{\xi})-e_{n}(\bm{x};\bm{\xi})}{\tau},v
\Big\rangle
+\Big\langle  e_{n+1} \frac{\partial e_{n+1}}{\partial \bm{x}}(\bm{x};\bm{\xi}),v
\Big\rangle\\
&\qquad\qquad+\left\langle \mathcal{L}(e_{n+1}(\bm{x};\bm{\xi}),u_{k-1,n+1}(\bm{x};\bm{\xi})),v\right\rangle 
=\langle r_{k,n+1}(\bm{x};\bm{\xi}),v\rangle.
\end{aligned}
$$
To further control computational complexity, we employ semi-implicit schemes \cite{L2006Semi‐implicit}.  Specifically, we use $e_{n}\frac{\partial e_{n+1}}{\partial \bm{x}}$ to approximate the nonlinear term $ e_{n+1}\frac{\partial e_{n+1}}{\partial \bm{x}}$ in the above equation.

Now we focus on obtaining the parameter-dependent basis function $\zeta_k(t;\bm{\xi})$ from \eqref{eq-frame-rode} with initial condition \eqref{eq-intial-rode}, i.e.,
\begin{equation}
\label{eq-nonlinear-xi}
\begin{aligned}
\Big\langle \frac{\partial (g_k(\bm{x},t)\zeta_k(t;\bm{\xi}))}{\partial t},g_k(\bm{x},t)\Big\rangle
+\Big\langle g_k(\bm{x},t)\zeta_k(t;\bm{\xi})\frac{\partial (g_k(\bm{x},t)\zeta_k(t;\bm{\xi}))}{\partial \bm{x}},g_k(\bm{x},t)\Big\rangle\\
+\left\langle\mathcal{L}(g_k(\bm{x},t)\zeta_k(t;\bm{\xi}),u_{k-1}(\bm{x},t;\bm{\xi})),g_k(\bm{x},t)\right\rangle
=\langle r_k(\bm{x},t;\bm{\xi}),g_k(\bm{x},t)\rangle,
\end{aligned}
\end{equation}
which is also a reformulation of \eqref{eq-frame-1} at step $k=1$ with $r_1(\bm{x},t;\bm{\xi})\equiv 0,$
and 
$$
\mathcal{L}(g_k(\bm{x},t)\zeta_k(t;\bm{\xi}),u_{k-1}(\bm{x},t;\bm{\xi}))\equiv-\mathcal{A}(g_1(\bm{x},t)\zeta_1(t;\bm{\xi});\bm{\xi}).
$$
Employing the backward Euler method for time discretization, \eqref{eq-nonlinear-xi} can be written as
$$
\begin{aligned}
    &\Big\langle \frac{g_{k,n+1}(\bm{x})-g_{k,n}(\bm{x})}{\tau}\zeta_{k,n+1}(\bm{\xi})+g_{k,n+1}(\bm{x})\frac{\zeta_{k,n+1}(\bm{\xi})-\zeta_{k,n}(\bm{\xi})}{\tau},g_{k,n+1}(\bm{x}) \Big\rangle
    \\[1mm]
    &+\Big\langle g_{k,n+1}(\bm{x})\zeta_{k,n+1}(\bm{\xi})\frac{\partial (g_{k,n+1}(\bm{x})\zeta_{k,n+1}(\bm{\xi}))}{\partial \bm{x}},g_{k,n+1}(\bm{x})\Big\rangle
    \\[1mm]
    &+\left\langle \mathcal{L}(g_{k,n+1}(\bm{x})\zeta_{k,n+1}(\bm{\xi}),u_{k-1,n+1}(\bm{x};\bm{\xi})),g_{k,n+1}(\bm{x})\right\rangle 
    =\langle r_{k,n+1}(\bm{x};\bm{\xi}),g_{k,n+1}(\bm{x})\rangle.
    \end{aligned}
$$
Applying the same approximation strategy to handle the nonlinear term, i.e., replacing
$
g_{k,n+1}(\bm{x})\zeta_{k,n+1}(\bm{\xi})\frac{\partial (g_{k,n+1}(\bm{x})\zeta_{k,n+1}(\bm{\xi}))}{\partial \bm{x}}$ with
$ g_{k,n}(\bm{x})\zeta_{k,n}(\bm{\xi})\frac{\partial (g_{k,n+1}(\bm{x})\zeta_{k,n+1}(\bm{\xi}))}{\partial \bm{x}}
$,
we can get a linear affine expression for $\zeta_{k,n+1}(\bm{\xi})$, which is omitted here as it is similar to \eqref{eq-iter-t} for the linear case.

Similarly to the linear case in \cref{sec-dvs-linear}, in practical computation, 
it is also not necessary to store all the values of $\zeta_{k}(t;\bm{\xi})$ at discrete times here but to use one simple ODE to get the value of $\zeta_{k}(t;\bm{\bar \xi})$ in the online stage,  which is based on a few scalars stored from the offline stage. 
Finally, we should mention that the strategy for handling the nonlinear terms described above ensures that the proposed DVS method also applies to other nonlinear parameter-dependent dynamical systems.

\section{Numerical experiments}
\label{sec-numerical examples}
In this section, we present four numerical examples and the corresponding numerical results to demonstrate the applicability and efficiency of the proposed DVS method.
In \cref{sec-num-linear}, we consider the reaction-diffusion equation and the two-dimensional heat equation to illustrate the performance of the proposed DVS method
for linear parameter-dependent dynamical systems.
Furthermore, we also apply the DVS method for nonlinear parameter-dependent dynamical systems exemplified by the Burgers equation and the Allen-Cahn equation in \cref{sec-num-nonlinear}.
 
To quantify the accuracy of the proposed DVS method, we use the average relative error $\epsilon$ for the solution defined by 
\begin{equation}
\label{eq-meanerror}
\epsilon:=\frac{1}{M}\sum_{i=1}^{M}\frac{||u(\bm{x},t;\bm{\xi}_i)-u_{N}(\bm{x},t;\bm{\xi}_i)||_{L^2([0,T];\mathcal{V})}}{||u(\bm{x},t;\bm{\xi}_i)||_{L^2([0,T];\mathcal{V})}},
\end{equation}
where $M$ is the number of samples used to compute the average error, 
$u_{N}(\bm{x},t;\bm{\xi})$ is the approximate solution obtained by the DVS method, and $u(\bm{x},t;\bm{\xi})$ is the reference solution.
In the following numerical examples, the reference solution is calculated using the finite element method in space and the backward Euler scheme in time (denoted FEM-BE).  

\subsection{Linear parameter-dependent dynamical systems}
\label{sec-num-linear}
In this subsection, we consider two linear dynamical systems: the reaction-diffusion equation with the parameter-dependent boundary condition and the heat equation with a parameter-dependent source term.

\subsubsection{Reaction-diffusion equation}
Consider the reaction-diffusion equation defined on domain $D=[0,1]$ given as follows
$$
\left\{
    \begin{aligned}
\frac{\partial u}{\partial t}+\xi_1u=&2\xi_2 \frac{\partial^2 u}{\partial {x}^2}+\xi_3, ~~ {x}\in D,~t\in[0,T],\\[.5mm]
u({x},0;\bm{\xi})=&2({x}+1)\xi_4, ~~ {x}\in D,\\[.5mm]
 u({x},t;\bm{\xi})=&2({x}+1)\xi_4,~~ {x}\in \partial D, ~t\in[0,T],
 \end{aligned}
 \right.
$$
in which we take $T=1$. 
Here, the parameter $\bm{\xi}:=({\xi}_1,{\xi}_2,{\xi}_{3},\xi_4)\in\mathbb{R}^{4} $, and each ${\xi}_i$ is in the interval $[1,3]$.
The reference solution is calculated using the finite element method in space with a mesh size $h=0.02$ and the backward Euler scheme in time with a step size $\tau=10^{-3}$. 
We take the sample size $|\Xi|=11$ for the offline stage in \Cref{algorithm-dvs}, where $|\Xi|$ denotes the cardinality of the set $\Xi$.

\begin{figure}[htbp]
    \centering
    \subfigure[Comparison of the iteration point]{
        \includegraphics[width=2.4in, height=1.8in]{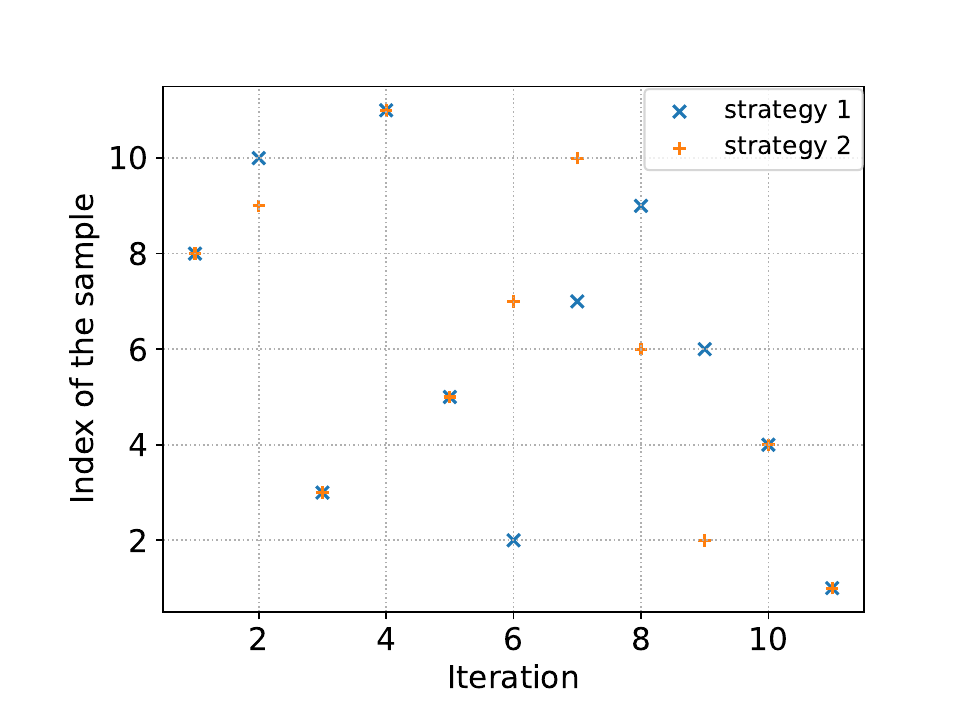}}
    \subfigure[Comparison of the average relative error]{
        \includegraphics[width=2.4in, height=1.8in]{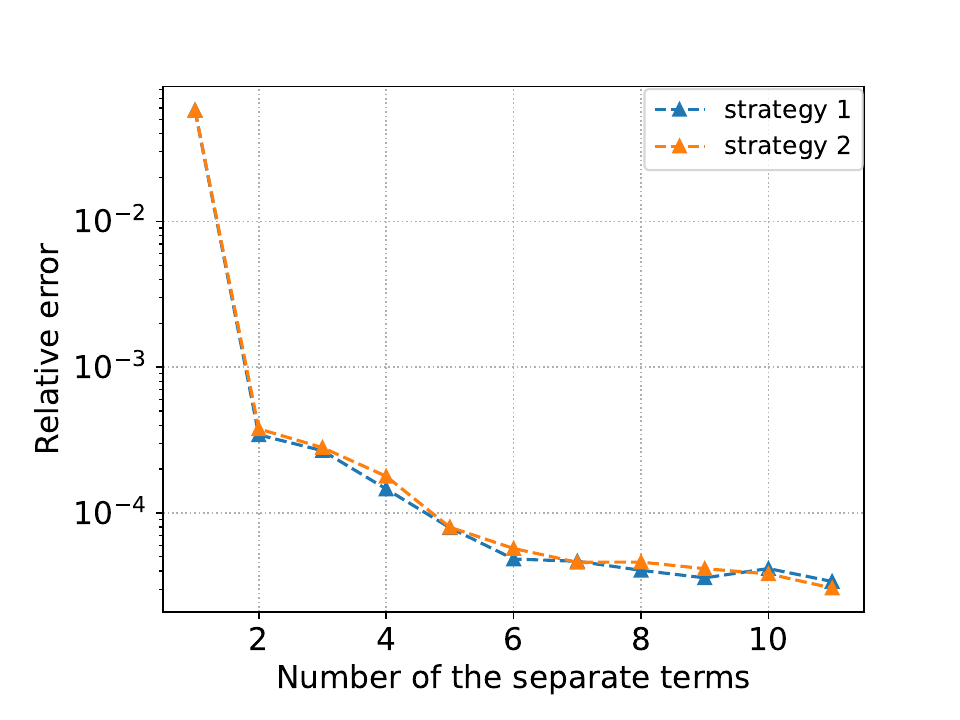}}
    \caption{Comparison of the two strategies to choose iteration points.}
    \label{fig-rd.1}
\end{figure}

First, we plot the comparison of the iteration point and the average relative error versus different numbers of separate terms $N$ in \Cref{fig-rd.1}. 
Here, the average relative error $\epsilon$ is calculated using equation (\ref{eq-meanerror}) with $M=10^3$ parameter samples. 

In \Cref{fig-rd.1}, we compared the two strategies for choosing the iteration point at each enrichment step in the offline stage, that is, by the approximate error $e(\bm{x},t;\bm{\xi})$ (strategy 1) and by the error estimator $\delta_k(t;\bm{\xi})$ (strategy 2) defined in \Cref{sec-error}. 
From the figure, we can see that (1) the two strategies make the same choice at certain iteration steps; 
(2) strategy 1 always gives a little better approximations than strategy 2, and both strategies have similar curves of error reduction, achieving the same level of approximation.  
In summary, using an error estimator to choose iteration points is appropriate in practical simulations. 

\begin{figure}[htbp]
    \centering
    \subfigure[Comparison of the relative error corresponding to the number of the separate terms being $N=2$, $N=4$, $N=6$, and $N=8$]{
     \includegraphics[width=2.4in, height=1.8in]{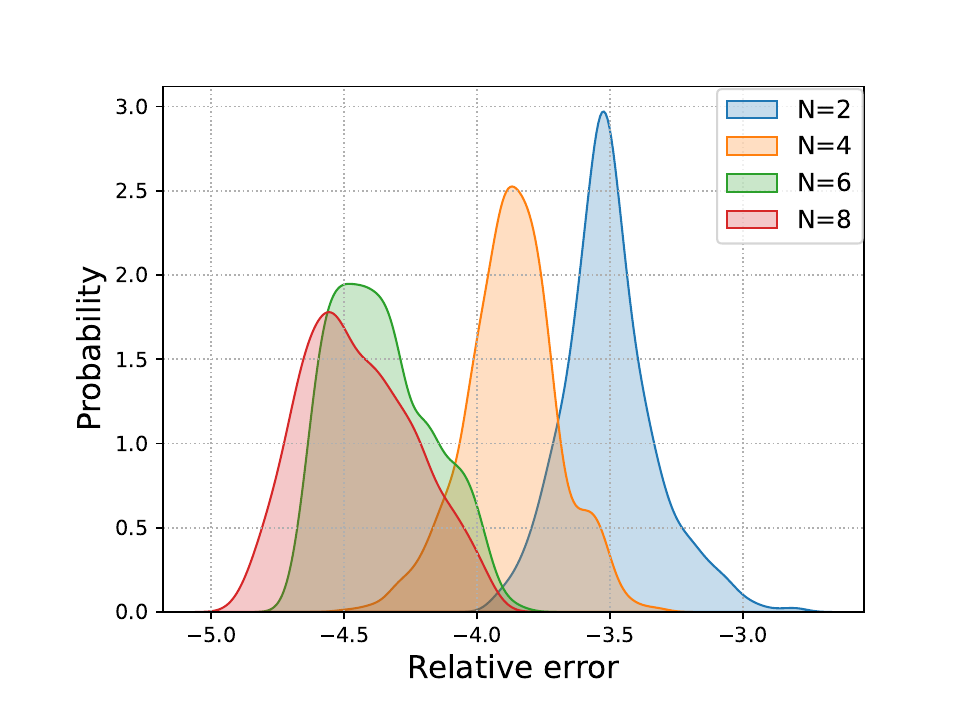}}
     \quad 
    \subfigure[The average relative error versus different time levels with the number of the separate terms being $N=1$, $N=3$, and $N=6$]{
        \includegraphics[width=2.4in, height=1.8in]{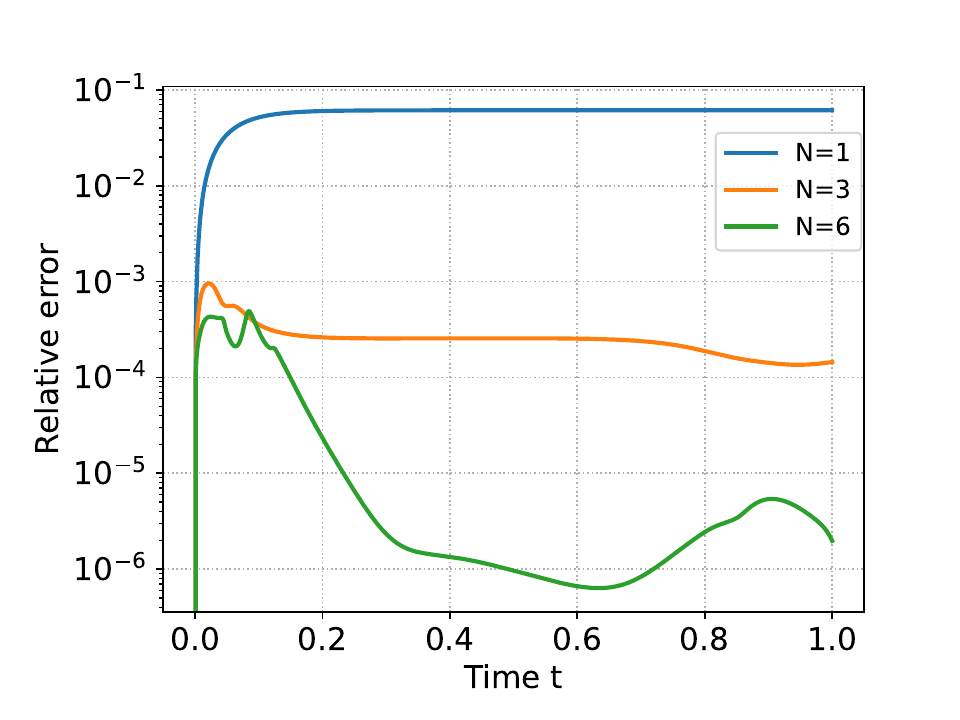}}
    \caption{Comparison of the relative error.}
    \label{fig-rd.2}
\end{figure}

To visualize individual relative error, we plot the probability density of logarithmic relative errors (with base 10) for each parameter sample with the number of the separate terms being  $N=2$, $N=4$, $N=6$, and $N=8$ in \Cref{fig-rd.2}(a). 
Here, the iteration point in the offline stage is chosen by the approximate error, and the following discussion is also based on this situation.
From this figure, we can see that the errors are relatively concentrated for a fixed $N$, and as the number of separate terms increases, the error decreases overall.
In \Cref{fig-rd.2}(b), we present the average relative error versus different time levels with the number of the separate terms being $N=1$, $N=3$, and $N=6$. 
From the figure, we know that as the number of the separate terms $N$ increases, the DVS method can provide better approximation at most fixed times.

\begin{figure}[htbp]
    \centering
    \subfigure[$g_1(\bm{x},t)$]{
        \includegraphics[width=1.825in, height=1.415in]{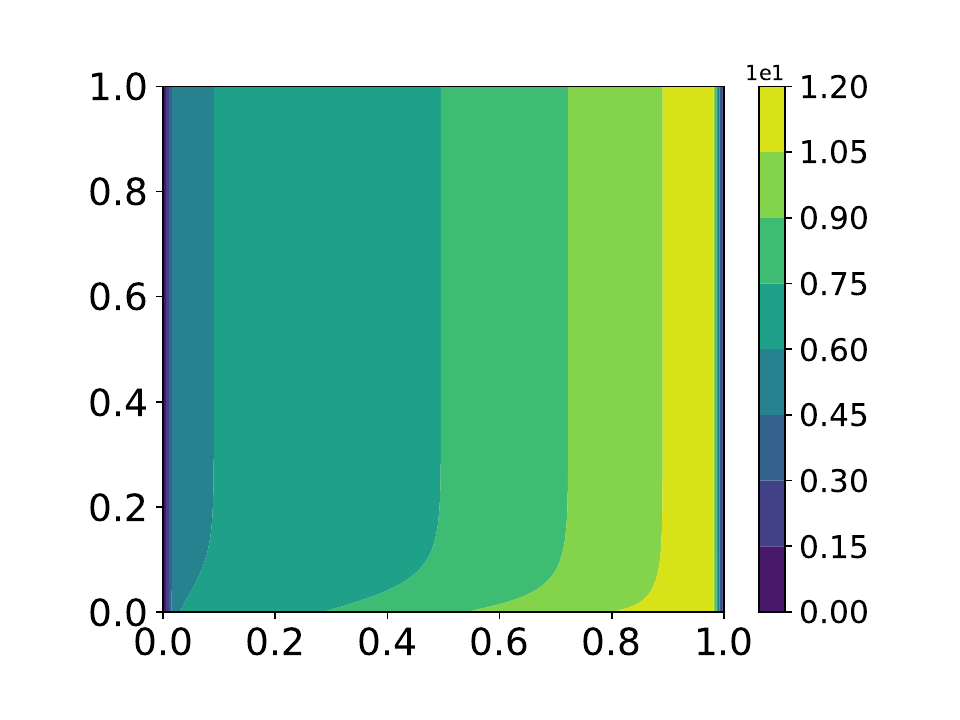}\hspace{-5mm}}
    \subfigure[$g_2(\bm{x},t)$]{
        \includegraphics[width=1.825in, height=1.415in]{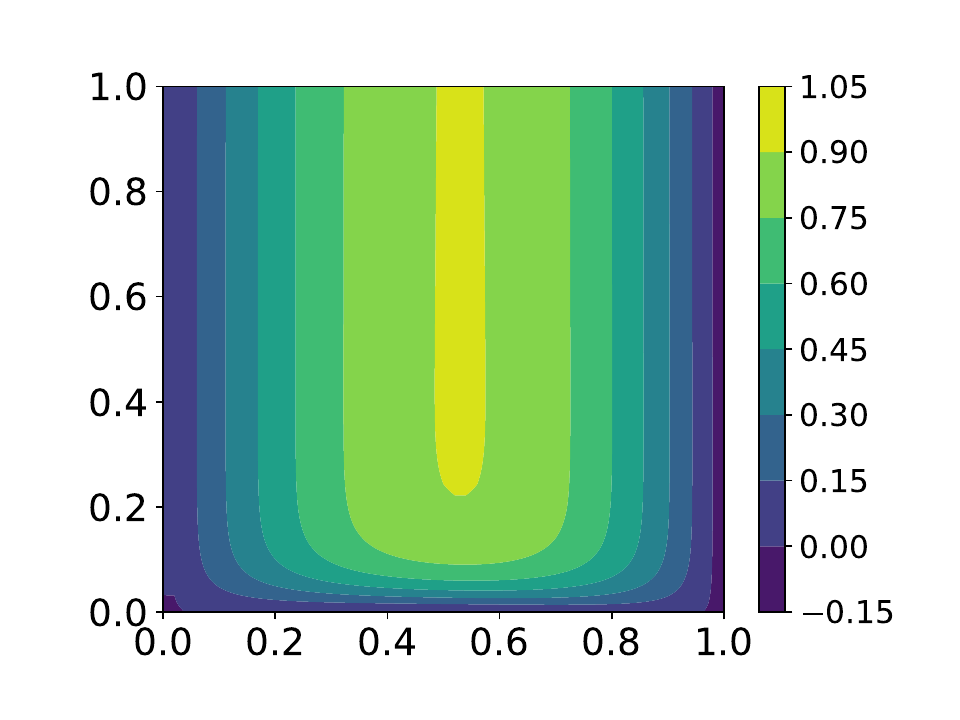}\hspace{-5mm}}
    \subfigure[$g_3(\bm{x},t)$]{
        \includegraphics[width=1.825in, height=1.415in]{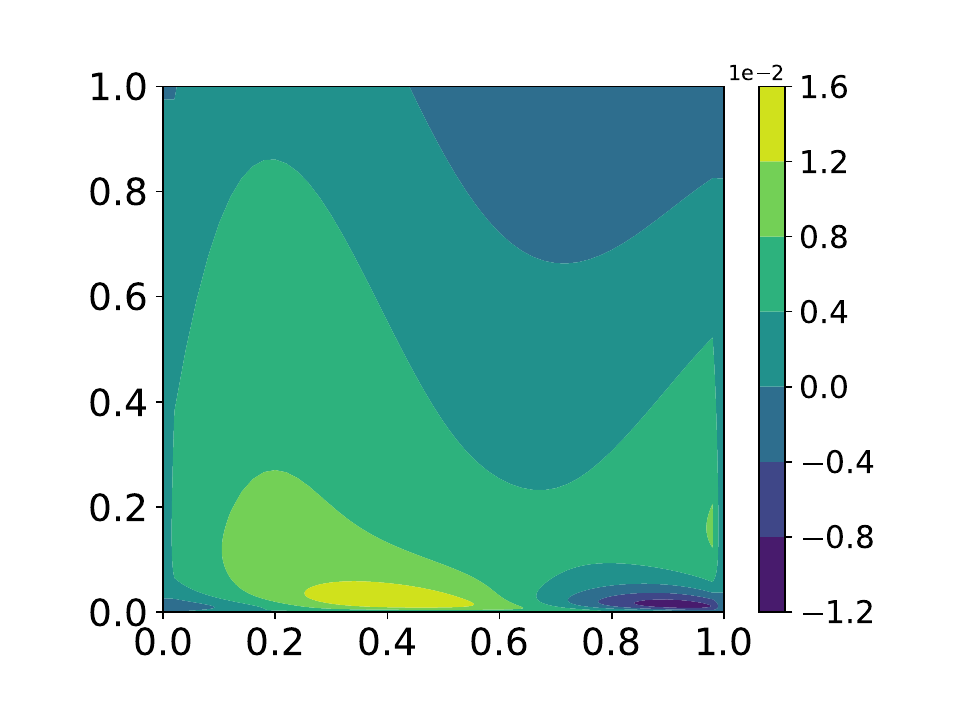}\hspace{-5mm}}\\
    \subfigure[$g_5(\bm{x},t)$]{
        \includegraphics[width=1.825in, height=1.415in]{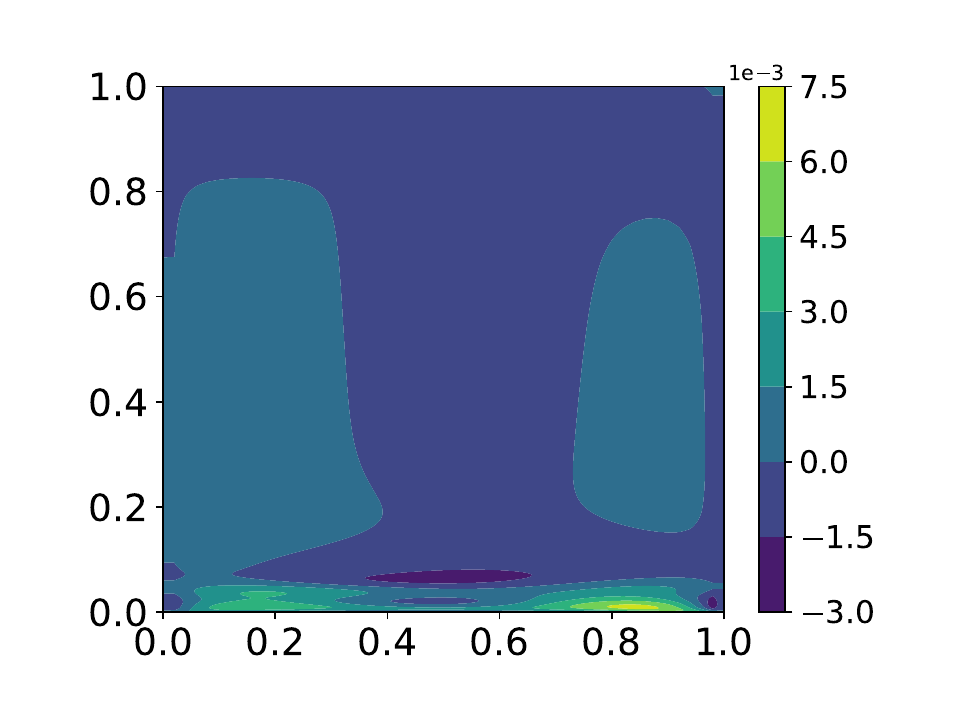}\hspace{-5mm}}
    \subfigure[$g_7(\bm{x},t)$]{
        \includegraphics[width=1.825in, height=1.415in]{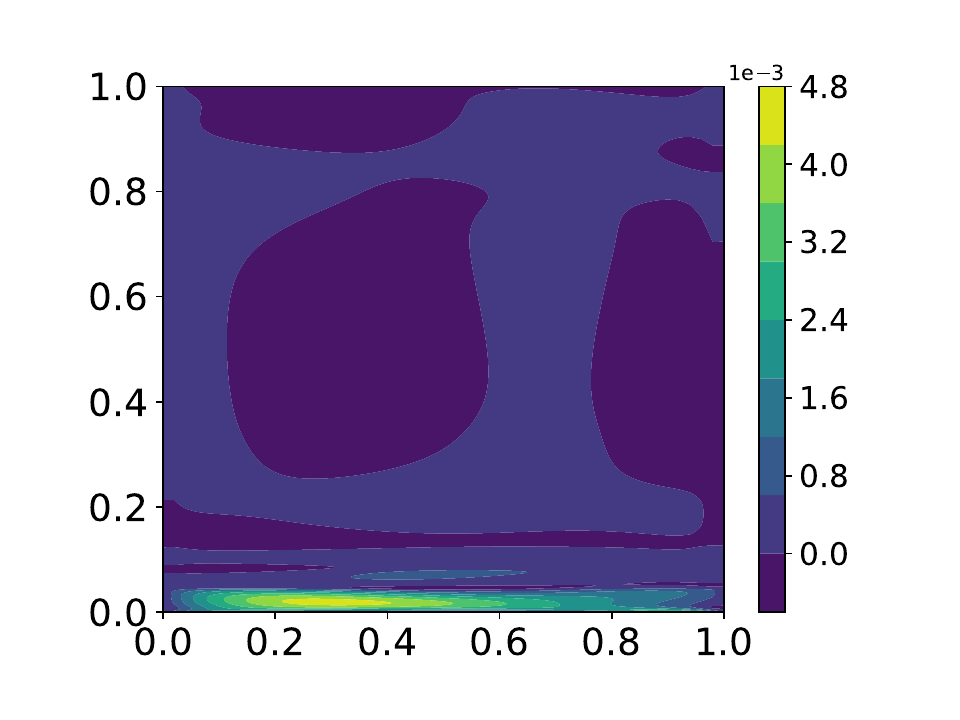}\hspace{-5mm}}
    \subfigure[$g_9(\bm{x},t)$]{
        \includegraphics[width=1.825in, height=1.415in]{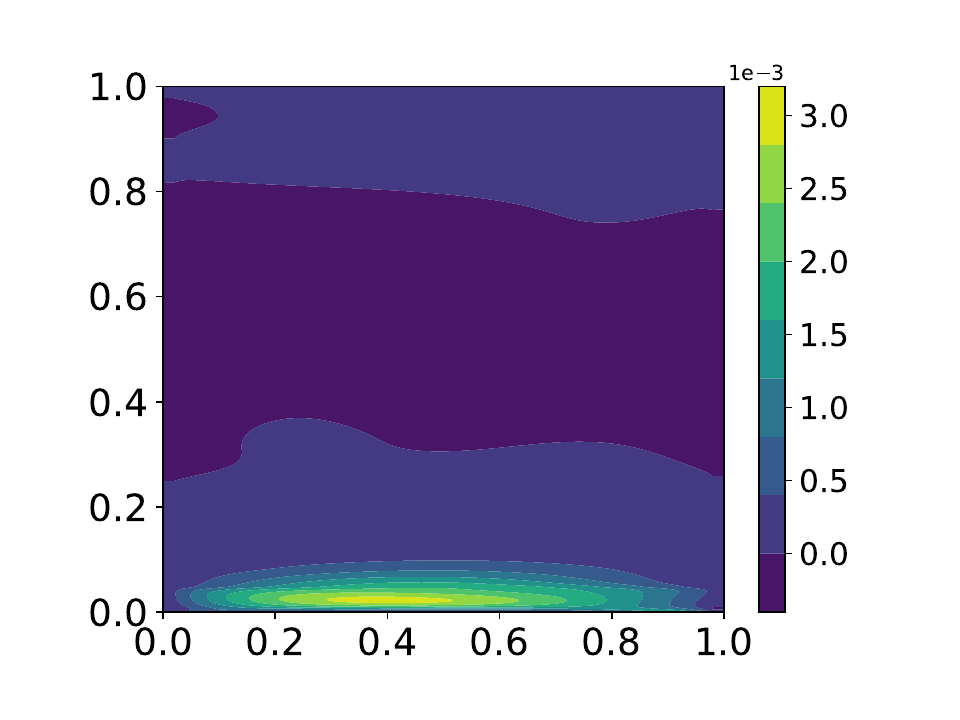}\hspace{-5mm}}
    \caption{Basis fields $g_i(\bm{x},t)$ of the DVS approximation.}
    \label{fig_rd.4}
\end{figure}

The basis fields $g_i(\bm{x},t)$ are shown in \Cref{fig_rd.4}, which indicates that the first field $g_1(\bm{x},t)$ represents the core information of the solution $u(\bm{x},t;\bm{\xi})$, and the last few fields capture the fine-scale information of the original system.

\begin{figure}[htbp]
    \centering
    \subfigure[ $N=1,2, \cdots, 10$]{
        \includegraphics[width=2.4in, height=1.8in]{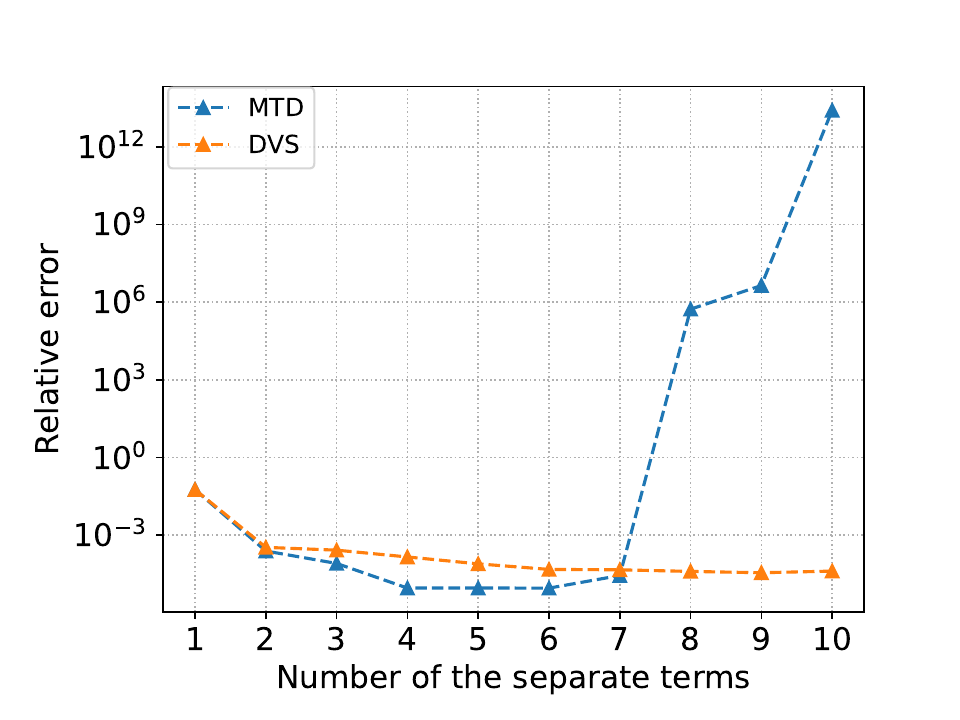}}
    \subfigure[ $N=1,2, \cdots, 7$]{
        \includegraphics[width=2.4in, height=1.8in]{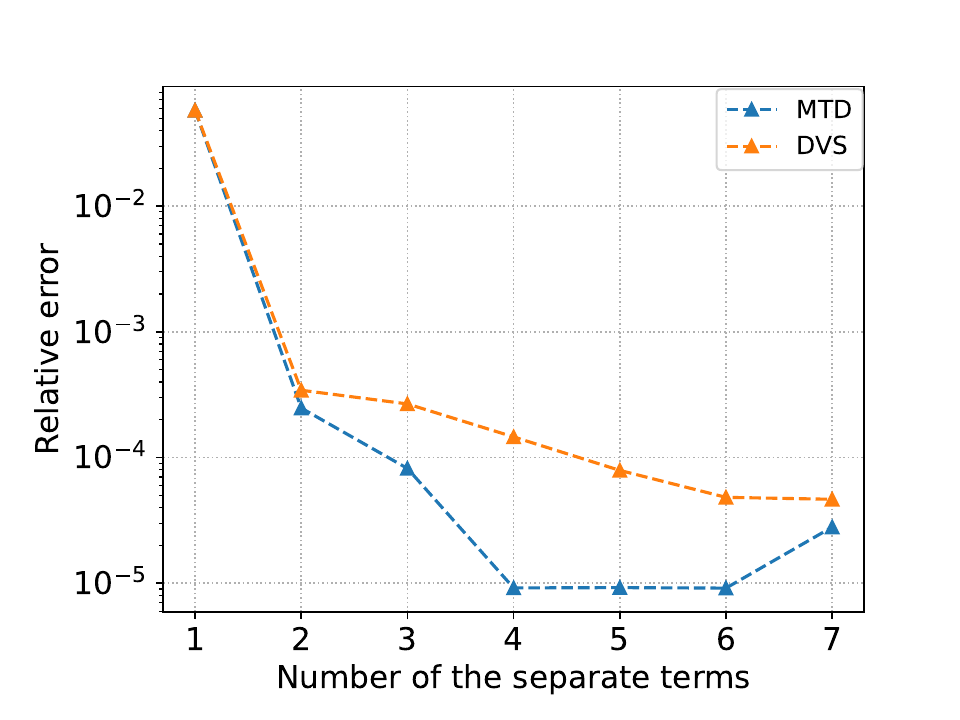}}
    \caption{The average relative error corresponding to the number of the separate terms $N$.}
    \label{fig_MTD}
\end{figure}

Next, we give a comparison of the DVS method with the model order reduction method using a time-dependent reduced basis (abbreviated as MTD) in \cite{Friess2017DynamicalMR}. 
We first plot the average relative error versus different numbers of separate terms $N$ for DVS and MTD in \Cref{fig_MTD}.
The figure on the left is the average relative error versus $N=1,2, \cdots, 10$, and the figure on the right is the average relative error versus $N=1,2, \cdots, 7$.
From the figure, we can see that:
(1) after $N=7$, the relative error of the MTD method increases quickly,
which may be due to the instability of the linear systems of equations to obtain the parametric coefficients (corresponding to the parametric basis $\zeta_i(t;\bm{\xi})$ in DVS), whose coefficient matrices are likely to have large condition numbers;
(2) the relative error of MTD decreases faster than that of DVS. 
This is probably because all parametric coefficients in MTD are updated in each iteration, while DVS does not modify the time-parameter basis function computed in the previous iterations.  
In \Cref{table.1}, we list the average relative error along with the CPU time based on $M=10^3$ parameter samples. The CPU time includes the offline CPU time ($\mathcal{T}_{\text{off}}$), online CPU time ($\mathcal{T}_{\text{on}}$), total CPU time ($\mathcal{T}_{\text{tot}}$), and the average online CPU time ($T_{\text{on}}$). 
From the table, we know that:
(1) the relative error with the same separate terms $N$ obtained by MTD is less than that of DVS;
(2) the magnitude of the average online CPU time for MTD and DVS is significantly smaller than that of FEM-BE;
(3) DVS needs less CPU time than MTD for the same $N$.
In summary, MTD provides higher accuracy under certain conditions, and DVS performs better in computational efficiency and error stability, making it a compelling choice for solving parameter-dependent dynamical systems.
Similar phenomena also occur in the remaining numerical experiments, so we omit to report these results and show only the effectiveness of DVS.

\begin{table}[h]
	\centering
	\caption{Comparison of the average relative errors and CPU time for MTD, DVS, and FEM-BE with different numbers of separate terms.}
        \vspace*{2pt}
        \scriptsize
        \begin{tabular}{cc|c|c|c|c|c}
			\Xhline{1pt}
			\multicolumn{2}{c|}{\centering Algorithm}
			&
			\multicolumn{1}{c|}{\centering $\varepsilon$}
			&
			\multicolumn{1}{c|}{\centering          $\mathcal{T}_{\text{off}}$}
                &
			\multicolumn{1}{c|}{\centering $\mathcal{T}_{\text{on}}$}
                &
			\multicolumn{1}{c|}{{\centering $\mathcal{T}_{\text{tot}}$}}
            &
			\multicolumn{1}{c}{{$T_{\text{on}}$ }}
                \\ 
			\cline{1-7}	
            \multirow{2}{*}{$N=2$} 
            &
            {\centering MTD}
          &   $2.48\times10^{-4}$    &  $8.43 s$     &   $19.93s$    &   $28.36 s$&  $1.99\times10^{-2}s$\\
          \cline{3-7}
      &{\centering DVS}&  $3.43\times10^{-4}$   &   $8.01 s$ &  $1.98s$   &   $9.99 s$  &  $1.98\times10^{-3}s$\\
      \cline{1-7}
             \multirow{2}{*}{$N=4$} 
            &
            {\centering MTD}
          &   $9.19\times10^{-6}$    &  $18.78 s$     &   $22.48s$    &   $41.26 s$&  $2.25\times10^{-2}s$\\
          \cline{3-7}
      &{\centering DVS}&  $1.46\times10^{-4}$   &   $16.51 s$ &  $3.96s$   &   $20.47 s$  &  $3.96\times10^{-3}s$\\  
      \cline{1-7}
           \multirow{2}{*}{$N=7$} 
            &
            {\centering MTD}
          &   $2.80\times10^{-5}$    &  $39.92 s$     &   $118.52s$    &   $158.44 s$&  $1.19\times10^{-1}s$\\
          \cline{3-7}
      &{\centering DVS}&  $4.66\times10^{-5}$   &   $29.86 s$ &  $7.20s$   &   $37.06s$  &  $7.20\times10^{-3}s$  \\
      \cline{1-7}
    \multicolumn{2}{c|}{\centering FEM-BE} & $\setminus $ & $\setminus $ & $\setminus $ &$3752.14s $&$3.75s $\\
	\Xhline{1pt}
		\end{tabular}
	\label{table.1}
    \end{table}

\subsubsection{Heat equation}
\label{num-heat}
In this numerical experiment, we consider the heat equation defined on domain $D=[0,\pi]\times[0,\pi]$ given by
$$
    \left\{
    \begin{aligned}
        \frac{\partial u}{\partial t}&=\kappa(\bm{\xi})\Delta u+f, ~~ \bm{x}\in D,~t\in[0,T],\\
        u(\bm{x},0;\bm{\xi})&=\sin(x_1)\sin(x_2)+1, ~~ \bm{x}\in D,
        \\[1mm]
        u(\bm{x},t;\bm{\xi})&=1,~~ \bm{x}\in \partial D, ~t\in[0,T],
    \end{aligned}
    \right.
$$
where $\kappa(\bm{\xi})={\xi}_1$, $T=1$, and the source function is defined by 
$$
f(\bm{x};\bm{\xi})=1+\sum_{m=1}^{10}\frac{\sin(2\pi mx_1)+\sin(2\pi mx_2)}{m^2\pi^2}{\xi}_{m+1}.
$$
 Here, the parameter $\bm{\xi}:=({\xi}_1,{\xi}_2,\dots,{\xi}_{11})\in\mathbb{R}^{11} $, and each ${\xi}_i$, $i=1,2,\dots,11$, is in the interval $[1,4]$.
The reference solution is calculated using the finite element method in space with a mesh size $h_{{x}_1}=h_{{x}_2}=\pi/50$ and the backward Euler scheme in time with a step size $\tau=10^{-4}$. We choose $|\Xi|=12$ as the sample size for the offline stage in \Cref{algorithm-dvs}.

In order to illustrate the necessary of the time and parameter-dependent basis, we apply the VS method (instead of DVS) directly to construct the separate approximation such that  
\begin{equation}
    \label{eq-appro-heat2}
    u(\bm{x},t;\bm{\xi})\approx u_{N}(\bm{x},t;\bm{\xi}):=\sum_{i=1}^{N}\zeta_i(\bm{\xi})g_i(\bm{x},t),
\end{equation}
where each $g_i(\bm{x},t)$ is a time-dependent spatial basis, but the parameter-dependent basis $\zeta_i(\bm{\xi})$ is time-independent. This method belongs to linear model reduction methods, and the procedure is similar to the DVS method, more details can be found in \Cref{sec-vs}.

\begin{figure}[htbp]
    \centering
    \subfigure[DVS]{
        \includegraphics[width=2.4in, height=1.8in]{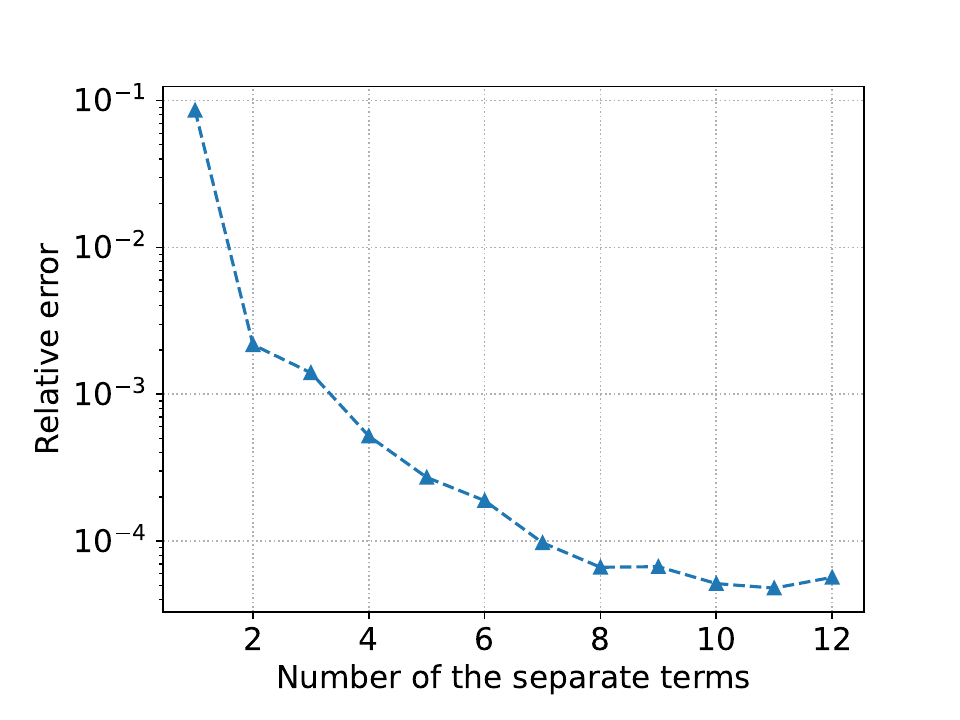}}
    \subfigure[VS]{
        \includegraphics[width=2.4in, height=1.8in]{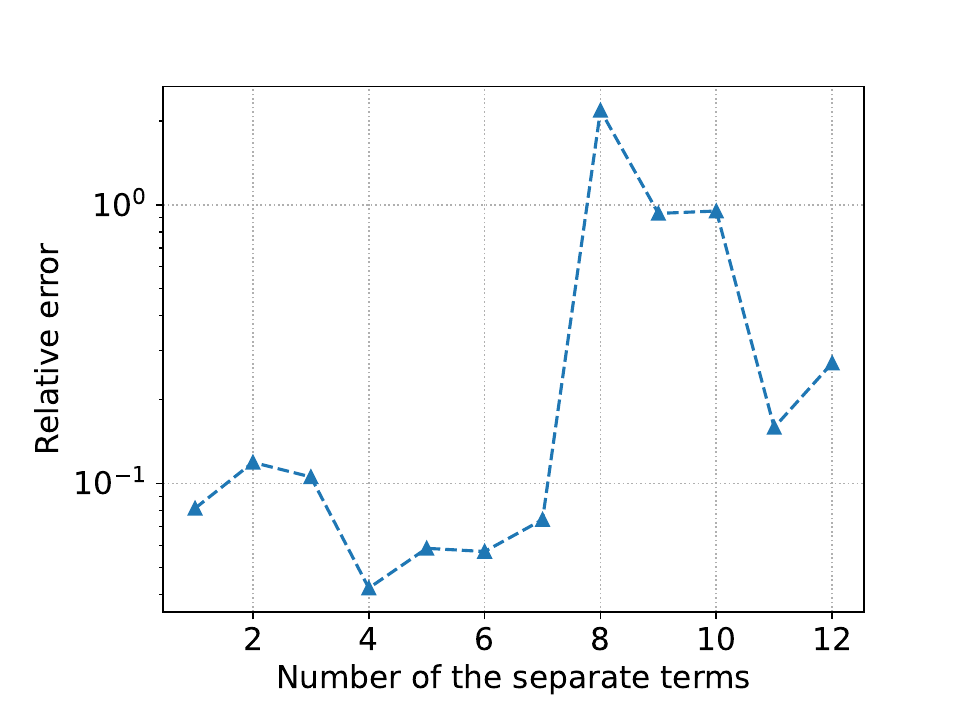}}
    \caption{The average relative error corresponding to the different numbers of the separate terms $N$.}
    \label{fig1.1}
\end{figure}

In \Cref{fig1.1}, we depict the average relative error of the approximate solution obtained by the DVS method and the VS method, that is, (\ref{eq-approx}) and (\ref{eq-appro-heat2}), versus different numbers of separate terms $N$, where the average relative error $\epsilon$ is calculated with $M=10^3$ parameter samples.  
The figures show that as the number of separate terms $N$ increases, the DVS method becomes more accurate. However, the average relative error of the VS method does not have this tendency. Furthermore, the DVS method can give a much better approximation accuracy than the VS method. 
	
\begin{figure}[htbp]
		\centering
         \subfigure[DVS]{
			\includegraphics[width=2.4in, height=1.8in]{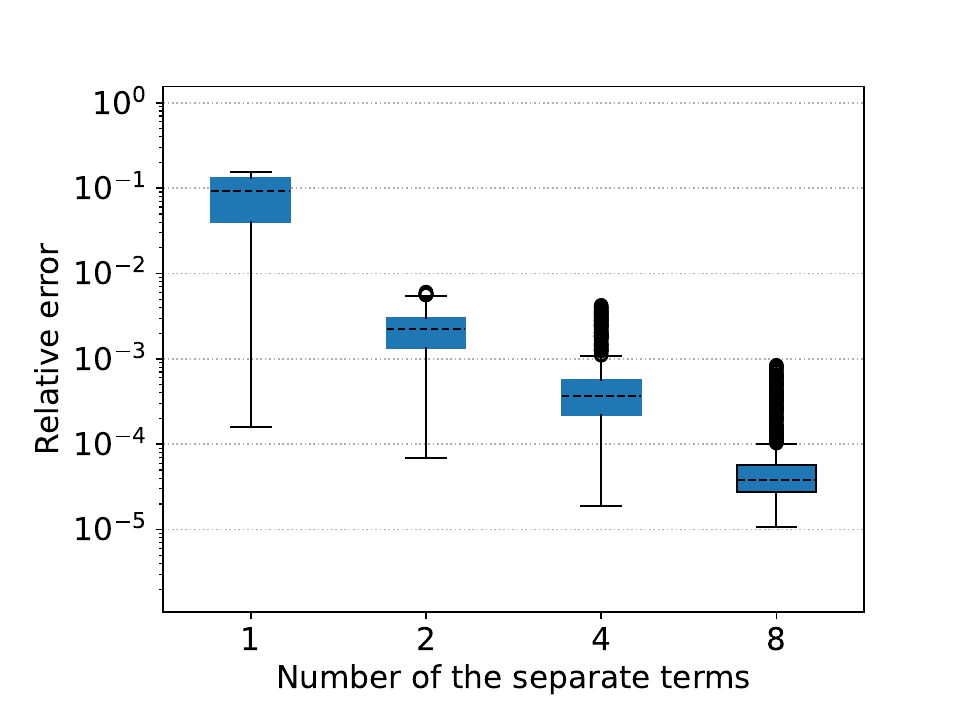}}
		\subfigure[VS]{
			\includegraphics[width=2.4in, height=1.8in]{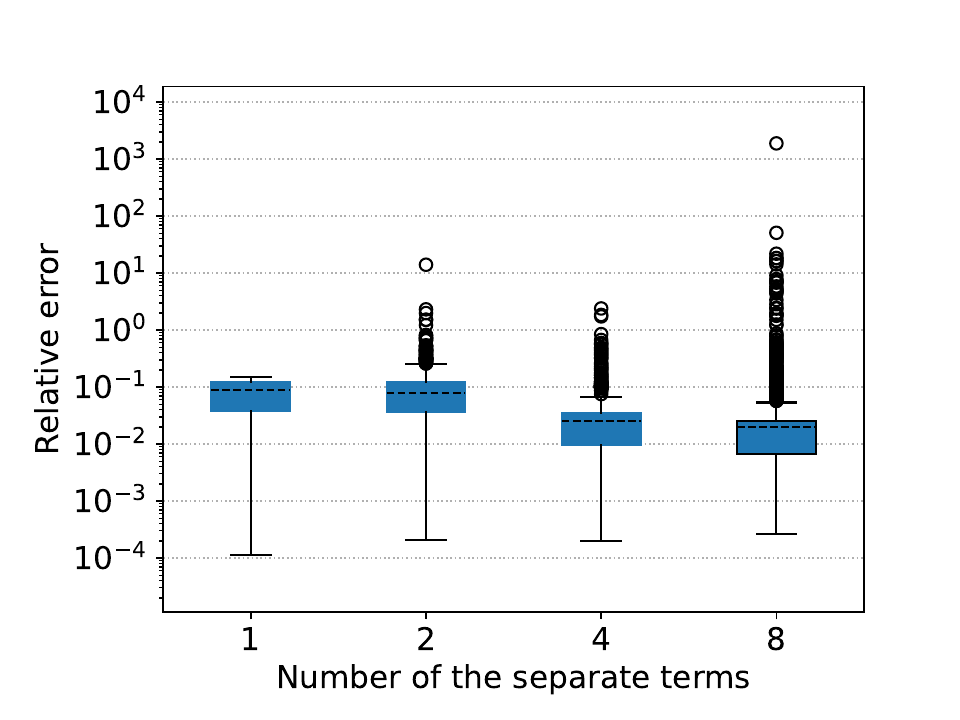}}
		\caption{The average relative errors with error bars corresponding to the number of the separate terms being $N=1$, $N=2$, $N=4$, and $N=8$.}
		\label{fig1.2}
	\end{figure}

    Next, we plot the average relative errors with error bars corresponding to the number of the separate terms being $N=1$, $N=2$, $N=4$, and $N=8$ in \Cref{fig1.2}. From the figures, we find that as the number of separate terms increases, the error bars of the relative errors for the DVS method are more compact than those of the VS method.
	
\begin{figure}[htbp]
		\centering
         \subfigure[DVS]{
			\includegraphics[width=2.4in, height=1.8in]{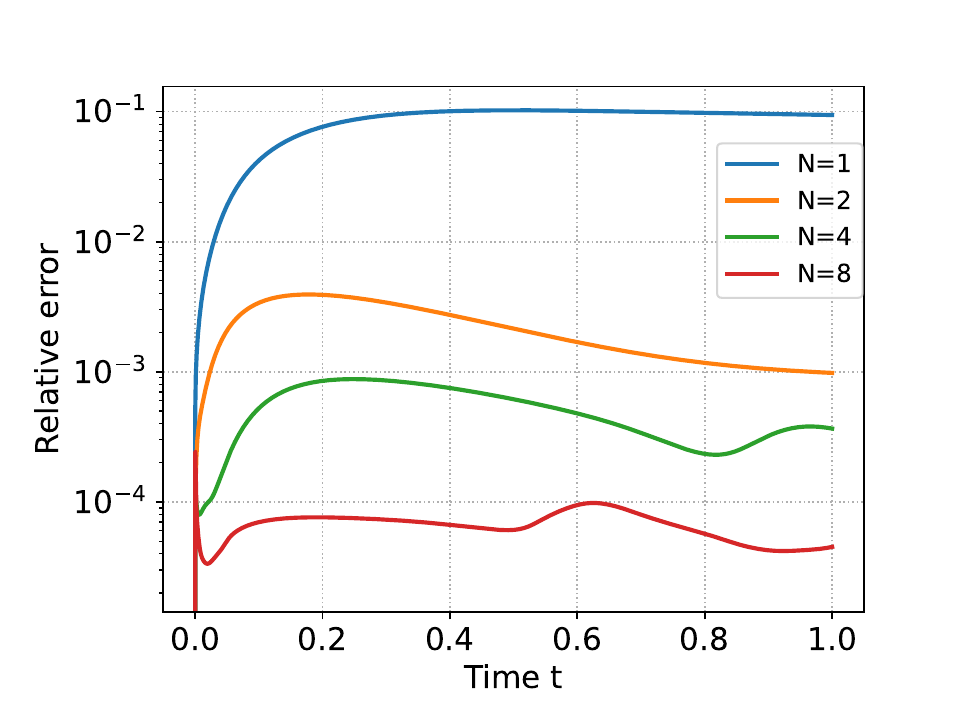}}
		\subfigure[VS]{
			\includegraphics[width=2.4in, height=1.8in]{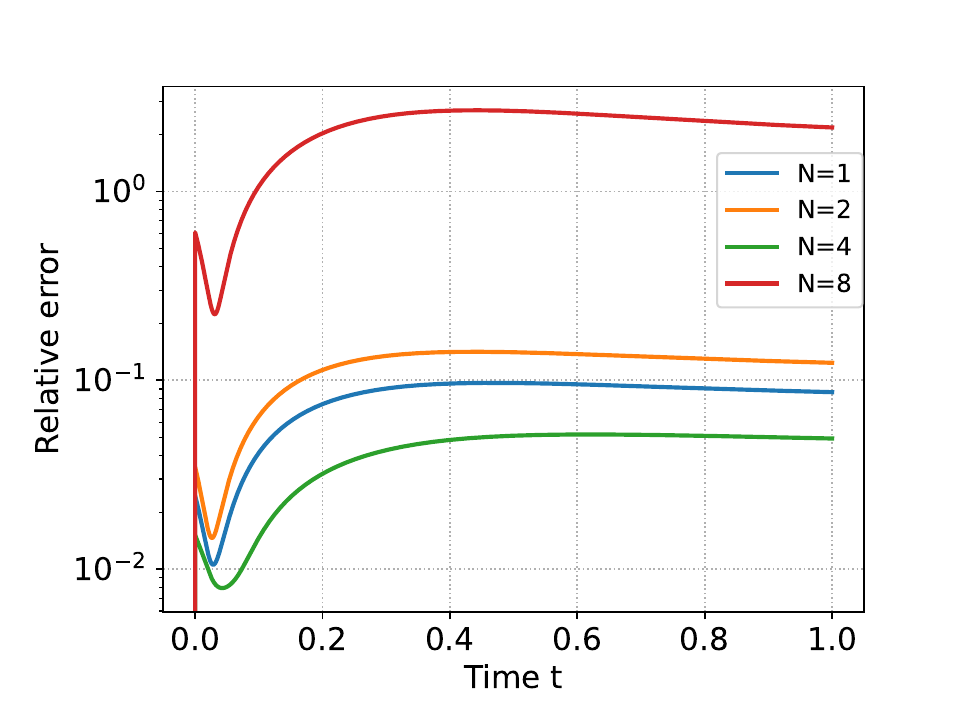}}
		\caption{The average relative error versus different time levels with the number of the separate terms being $N=1$, $N=2$, $N=4$, and $N=8$.}
		\label{fig1.3}
	\end{figure}

     \Cref{fig1.3} shows the average relative error versus different time levels with the number of separate terms being $N=1$, $N=2$, $N=4$, and $N=8$. From the figures, we see that at each fixed time $t$ the average relative error of the DVS method decreases more notably as the number of separate terms $N$ increases, compared to the VS method.

In \Cref{table_heat}, we list the average relative error and the average online CPU time based on $M=10^3$ parameter samples at fixed time $t=1$. From the table, we know that:
(1) the DVS method has a much better approximation than the VS method;
(2) the VS method needs less average online CPU time compared to the DVS method, and the magnitude of average online CPU time for both the DVS and VS methods is much smaller than that of the FEM-BE method;
(3) the DVS method achieves a much better trade-off in both approximation accuracy and computational efficiency than the VS method, which shows the necessity of using the time-dependent basis functions for both parameter and spatial variable.
	
    \begin{table}[hbtp]
	\centering
	\caption{Comparison of the average relative errors and online CPU time for DVS, VS and FEM-BE with different numbers of the separate terms.}
        \vspace*{2pt}
        \scriptsize
		\begin{tabular}{c|c|c|c|c}
			\Xhline{1pt}
			\multicolumn{1}{c|}{\multirow{2}{*}{\centering Algorithm }}
			&
			\multicolumn{2}{c|}{\centering DVS}
			&
			\multicolumn{2}{c}{\centering VS}\\
			\cline{2-5}	
			&
			\multicolumn{1}{c|}{\centering error $\epsilon$ }
			&
			\multicolumn{1}{c|}{\centering online time}
			&
			\multicolumn{1}{c|}{\centering error $\epsilon$ }
			&
			\multicolumn{1}{c}{\centering  online time}
                \\
			\hline
		   $N=2$    &  $9.81 \times10^{-4} $  &  $0.71s$  &  $1.24 \times10^{-1} $   &  $0.47s$ \\
			 $N=4$    &  $3.66 \times10^{-4} $  &  $2.37s$  &  $4.91 \times10^{-2} $   &  $1.52s$ \\
		   $N=6$    &  $1.77 \times10^{-4} $  &  $4.98s$  &  $6.37 \times10^{-2} $   &  $3.12s$\\
			 $N=8$    &  $4.52 \times10^{-5} $  &  $8.52s$  &  $2.18 $                 &  $5.28s$\\
			 $N=10$   &  $4.27 \times10^{-5} $  &  $13.01s$ &  $9.64 \times10^{-1} $   &  $8.00s$\\
			\hline
			\multicolumn{1}{c|}{\centering FEM-BE} & $\setminus $ &  $7.06\times10^{2}s$ &  $\setminus $   & $7.06\times10^{2}s$\\
			\Xhline{1pt}
		\end{tabular}
	\label{table_heat}
    \end{table}

\subsection{Nonlinear parameter-dependent dynamical systems}
    \label{sec-num-nonlinear}
    In this subsection, we will consider two nonlinear dynamical systems, including the Burgers equation with parameter-dependent initial conditions and the Allen-Cahn equation.
    
    \subsubsection{Burgers equation}
	Firstly, we consider the following Burgers equation
	$$
		\left\{
		\begin{aligned}
			\frac{\partial u}{\partial t}+u\frac{\partial u}{\partial {x}}&=\frac{{\xi}_1}{50}\frac{\partial^2 u}{\partial {x}^2},\ {x}\in [0,1],\ t\in[0,T],\\
			u({x},0;\bm{\xi})&=\frac{{x}(1-{x})}{2}{\xi}_2,\  {x}\in [0,1],\\
            u(0,t;\bm{\xi})&=u(1,t;\bm{\xi})=0,~~t\in[0,T],\\
		\end{aligned}
		\right.
	$$
	with $T=2$. The parameter $\bm{\xi}:=(\xi_1,\xi_2)\in \mathbb{R}^2$, and we set $\xi_i$, $i=1,2$ in the interval $[1,3]$.
	The reference solution is calculated using the finite element method in space with mesh size $h=0.01$ and the backward Euler scheme in time with step size $\tau=10^{-4}$, and the same strategy discussed in \cref{sec-dvs-burgers} is used to treat the nonlinear term.
    We choose $|\Xi|=12$ samples for training in \Cref{algorithm-dvs}.
	
\begin{figure}[htbp]
		\centering
		\includegraphics[width=.5\textwidth,height=1.8in]{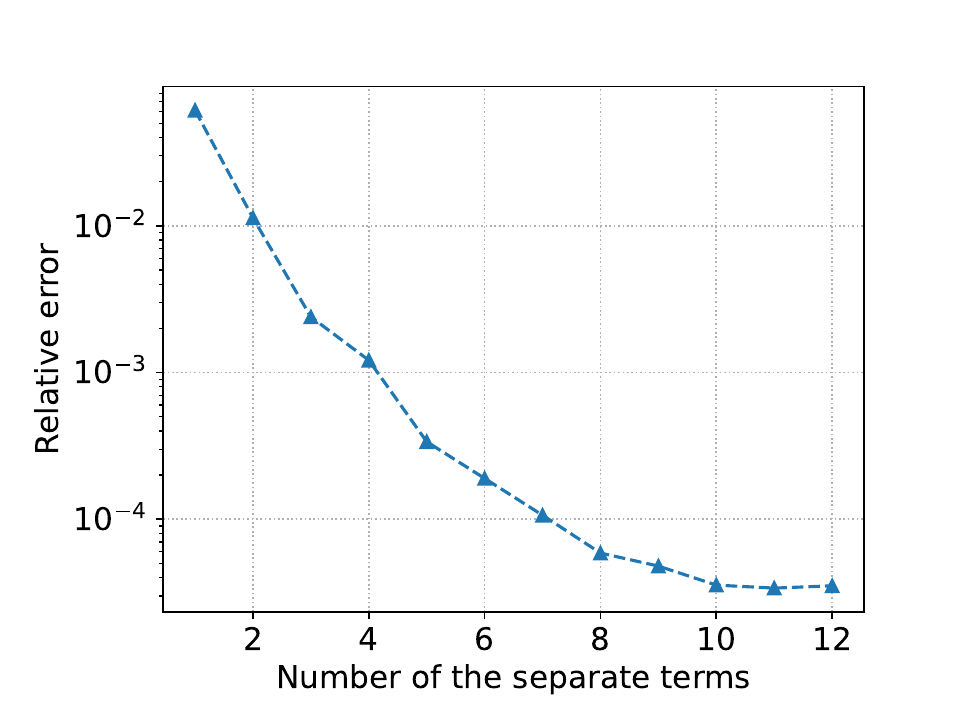}
		\caption{ The average relative error corresponding to the different numbers of the separate terms $N$.}
		\label{fig_burgers.1}
	\end{figure}

	Firstly, we represent the average relative error versus the number of separate terms $N$ for the DVS method in \Cref{fig_burgers.1}, where the average relative error $\epsilon$ is calculated based on $M=10^3$ test samples.
	From the figure, we know that as the number of separate terms $N$ increases, the approximation becomes more accurate.
	
\begin{figure}[htbp]
    \centering
    \subfigure[$t=1$ ]{
         \includegraphics[width=2.4in, height=1.8in]{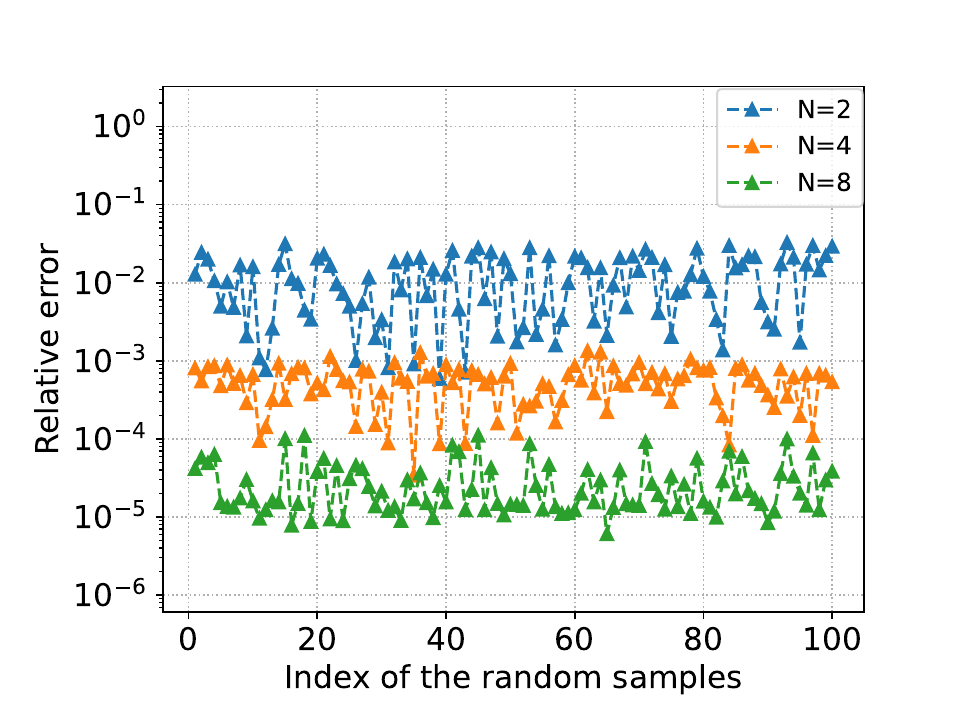}}
    \subfigure[$t=2$]{
        \includegraphics[width=2.4in, height=1.8in]{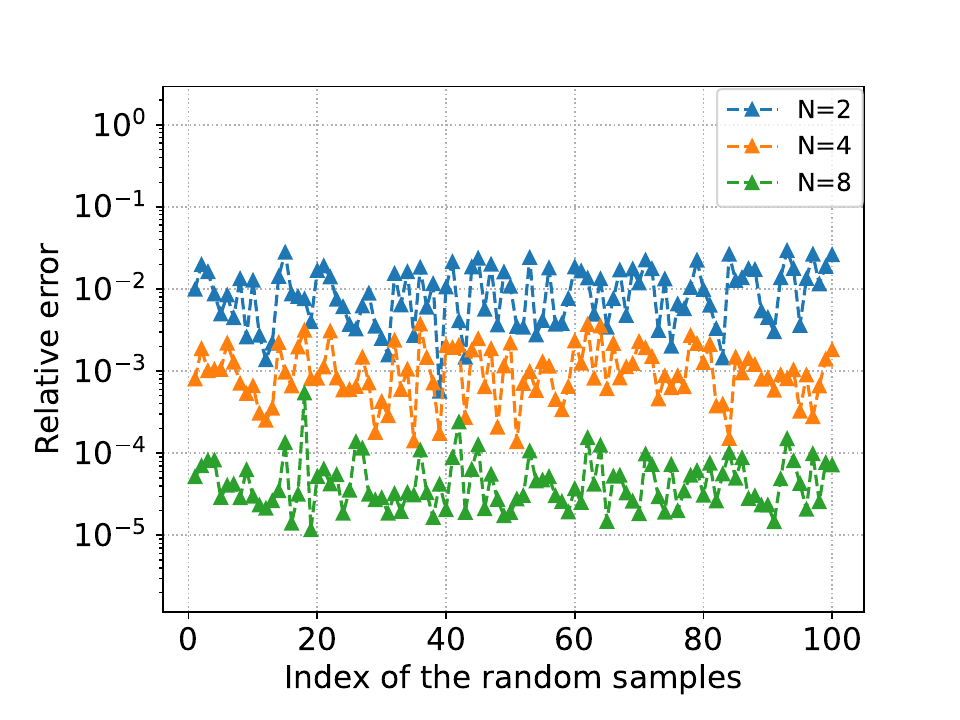}}
    \caption{The relative error for the first 100 samples versus the number of the separate terms being $N=2$, $N=4$, and $N=8$ at fixed time $t=1$ and $t=2$.}
    \label{fig_burgers.2}
\end{figure}

	In \Cref{fig_burgers.2}, we fix the number of separate terms as $N=2$, $N=4$, and $N=8$, and plot the relative error in fixed time $t=1$ and $t=2$ to visualize the individual relative error of the first $100$ samples from $M=10^3$ parameter samples. 
	From the figures, we see that the method we proposed gives a good approximation for each sample, and the individual relative error becomes smaller as the number of separate terms increases.

\begin{figure}[htbp]
		\centering
		\includegraphics[width=.5\textwidth,height=1.75in]{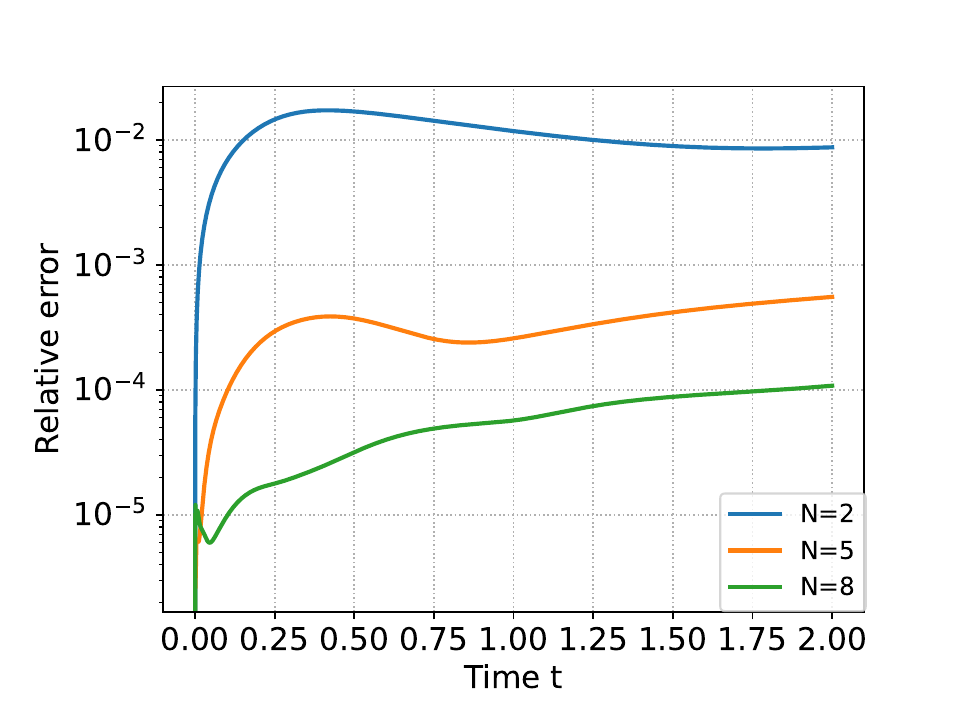}
		\caption{The average relative error for the solution versus different time levels with the number of the separate terms being $N=2$, $N=5$, and $N=8$.}
		\label{fig_burgers.3}
	\end{figure}

	The average relative error for the solution versus different time levels with the number of separate terms being $N=2$, $N=5$, and $N=8$ is shown in \Cref{fig_burgers.3}. This indicates that as the number of separate terms increases, the average relative error at fixed time becomes smaller.
	
	In \Cref{table_burgers}, we list the average relative error and the average online CPU time based on $M=10^3$ random samples at fixed times $t=1$ and $t=2$.
	From the table, we can conclude that (1) as the number of separate terms $N$ increases, the average online CPU time by the DVS method is added slowly; (2) the magnitude of average online CPU time by the DVS method is much smaller than that of the FEM-BE method;
(3) the approximation obtained by the DVS method can achieve a good trade-off in both approximation accuracy and computational efficiency.
\begin{table}[hbtp]
\centering
\caption{Comparison of the average relative errors and online CPU time for DVS and FEM-BE with different numbers of the separate terms.}
\vspace*{2pt}
\scriptsize
\begin{tabular}{cc|c|c|c|c}
\Xhline{1pt}
			\multicolumn{2}{c|}{\multirow{2}{*}{\centering Algorithm}}
			&
			\multicolumn{2}{c|}{\centering $t=1$ }
			&
			\multicolumn{2}{c}{\centering $t=2$}\\
		
			\cline{3-6}	
			&
			&
			\multicolumn{1}{c|}{\centering error $\epsilon$ }
			&
			\multicolumn{1}{c|}{\centering online time}
			&
			\multicolumn{1}{c|}{\centering error $\epsilon$ }
			&
			\multicolumn{1}{c}{\centering  online time}\\
			\hline
			\multirow{5}{*}{DVS}  & $N=2$    &  $1.17 \times10^{-2} $   &  $1.63 \times10^{-2}s$  &  $0.87 \times10^{-2} $   &  $2.92 \times10^{-2}s$\\
			\cline{3-6}
			                     & $N=4$    &  $1.09 \times10^{-3} $  &  $5.48 \times10^{-2}s$ &  $2.62 \times10^{-3} $   &  $9.87 \times10^{-2}s$ \\
			\cline{3-6}
			                     & $N=6$    &  $2.17 \times10^{-4} $  &  $1.23 \times10^{-1}s$  &  $4.10 \times10^{-4} $   &  $2.17 \times10^{-1}s$\\
			\cline{3-6}
			                     & $N=8$    &  $5.68 \times10^{-5} $  &  $2.18 \times10^{-1}s$  &  $1.07 \times10^{-4} $   &  $3.82 \times10^{-1}s$\\
			\cline{3-6}
			                     & $N=10$   &  $2.76 \times10^{-5} $  &  $3.44 \times10^{-1}s$ &  $7.12 \times10^{-5} $   &  $5.92 \times10^{-1}s$\\
			\hline
			\multicolumn{2}{c|}{\centering FEM-BE} & $\setminus $ &  $10.72s$ &  $\setminus $   &  $21.47s$\\
			\Xhline{1pt}
		\end{tabular}
		\label{table_burgers}
	\end{table}

\begin{figure}[htbp]
		\centering
		\subfigure[$g_1(\bm{x},t)$]{
			\includegraphics[width=1.825in, height=1.415in]{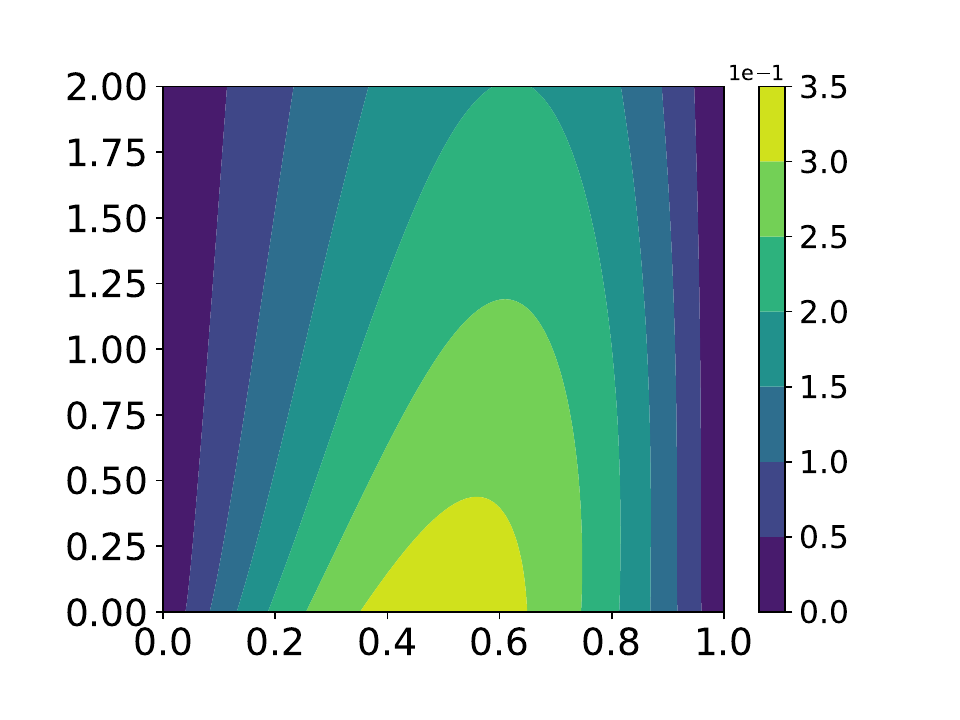}\hspace{-5mm}}
		\subfigure[$g_2(\bm{x},t)$]{
			\includegraphics[width=1.825in, height=1.415in]{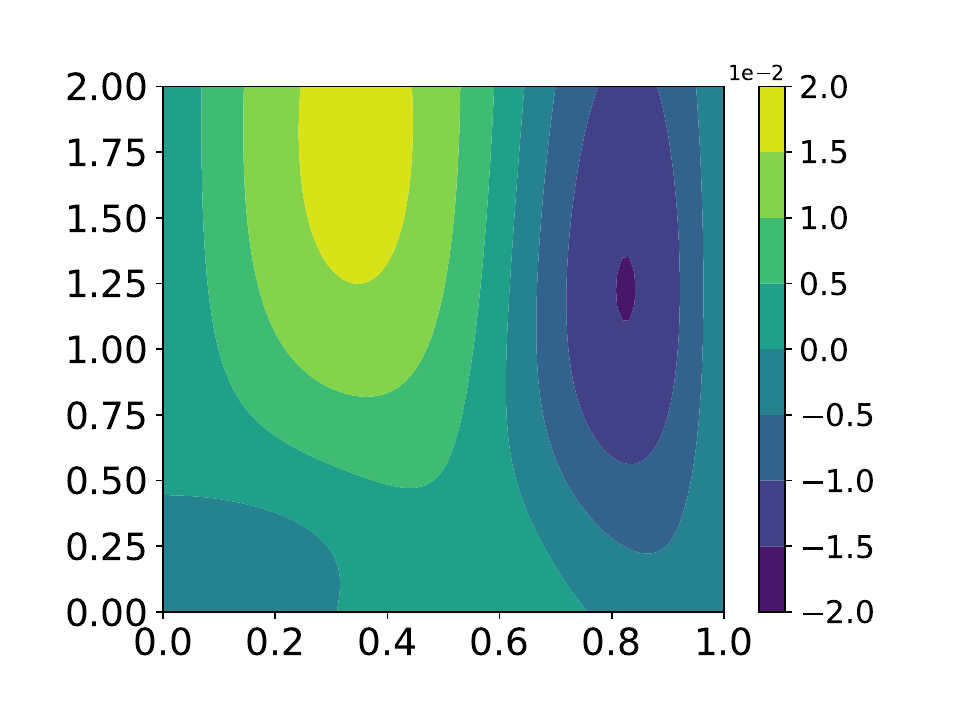}\hspace{-5mm}}
		\subfigure[$g_3(\bm{x},t)$]{
			\includegraphics[width=1.825in, height=1.415in]{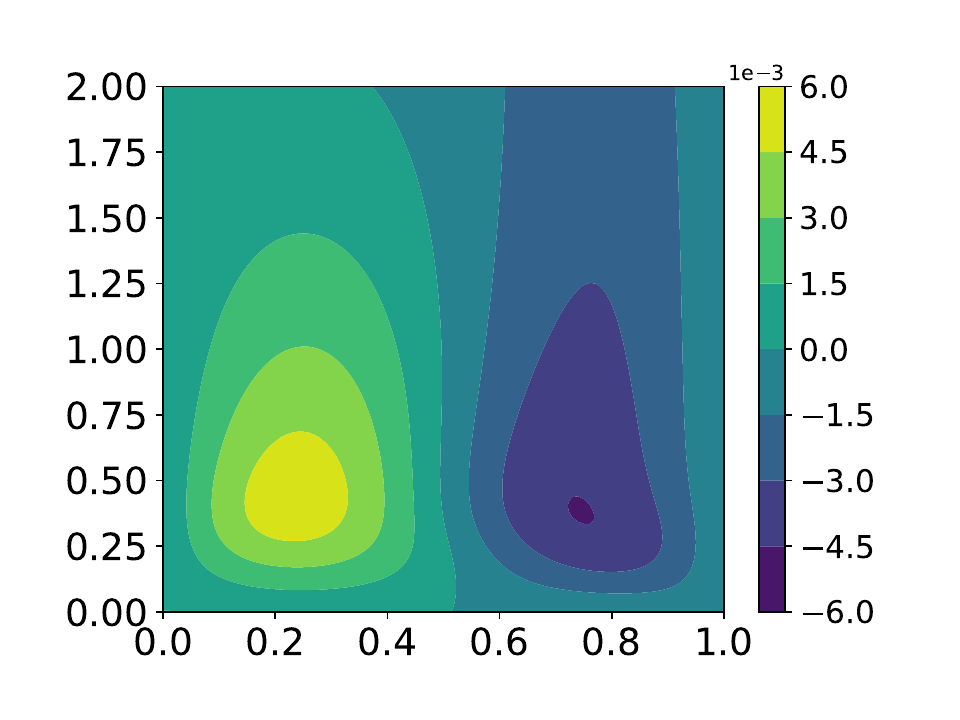}\hspace{-5mm}}\\
		\subfigure[$g_5(\bm{x},t)$]{
			\includegraphics[width=1.825in, height=1.415in]{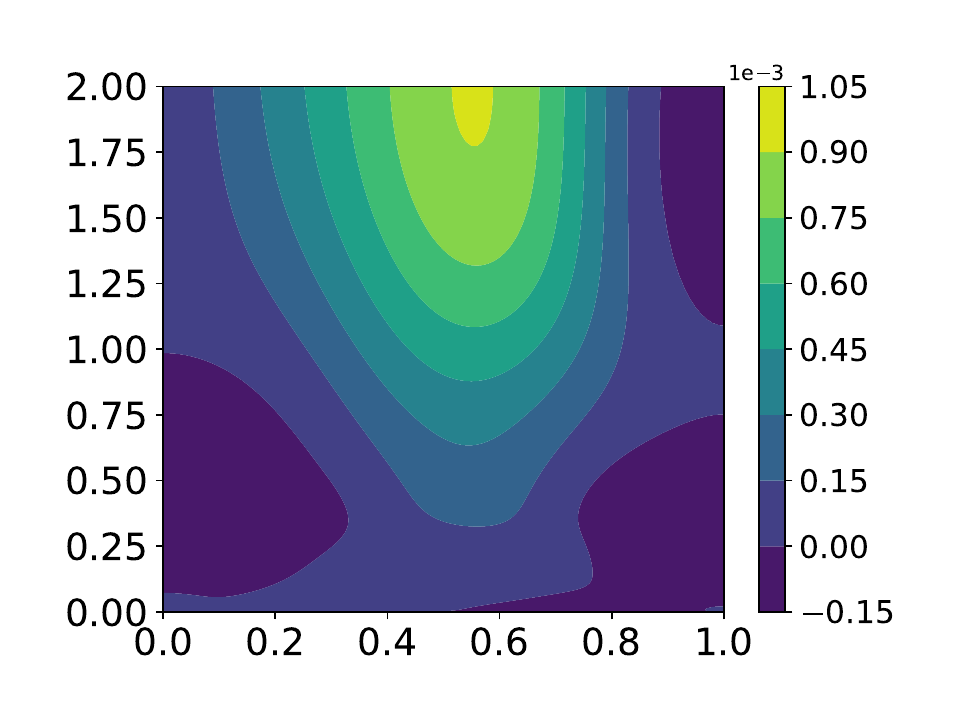}\hspace{-5mm}}
		\subfigure[$g_7(\bm{x},t)$]{
			\includegraphics[width=1.825in, height=1.415in]{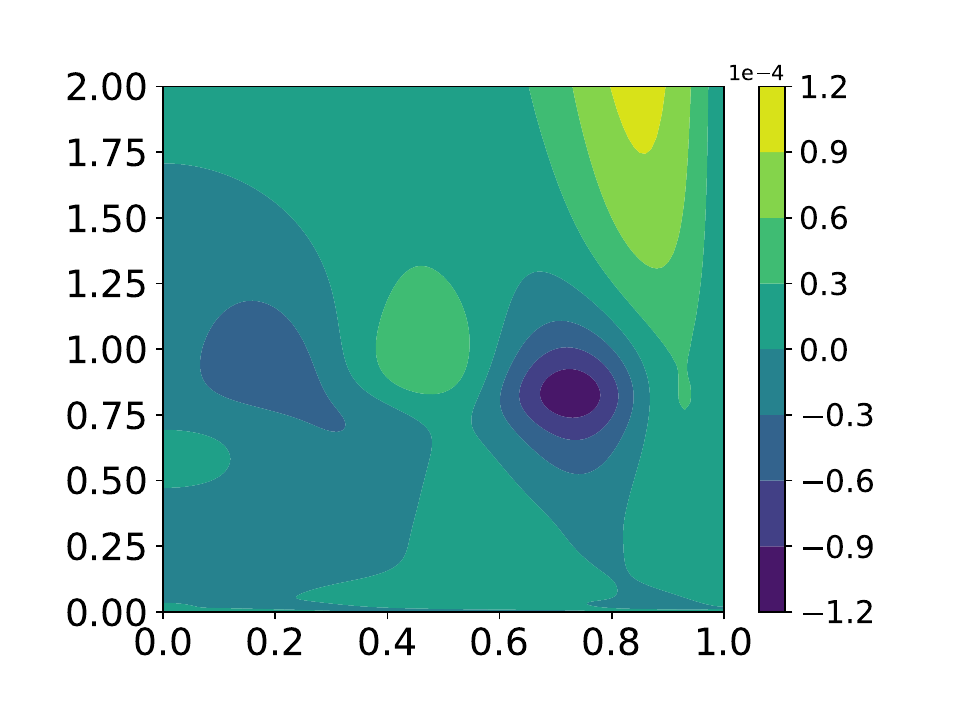}\hspace{-5mm}}
		\subfigure[$g_9(\bm{x},t)$]{
			\includegraphics[width=1.825in, height=1.415in]{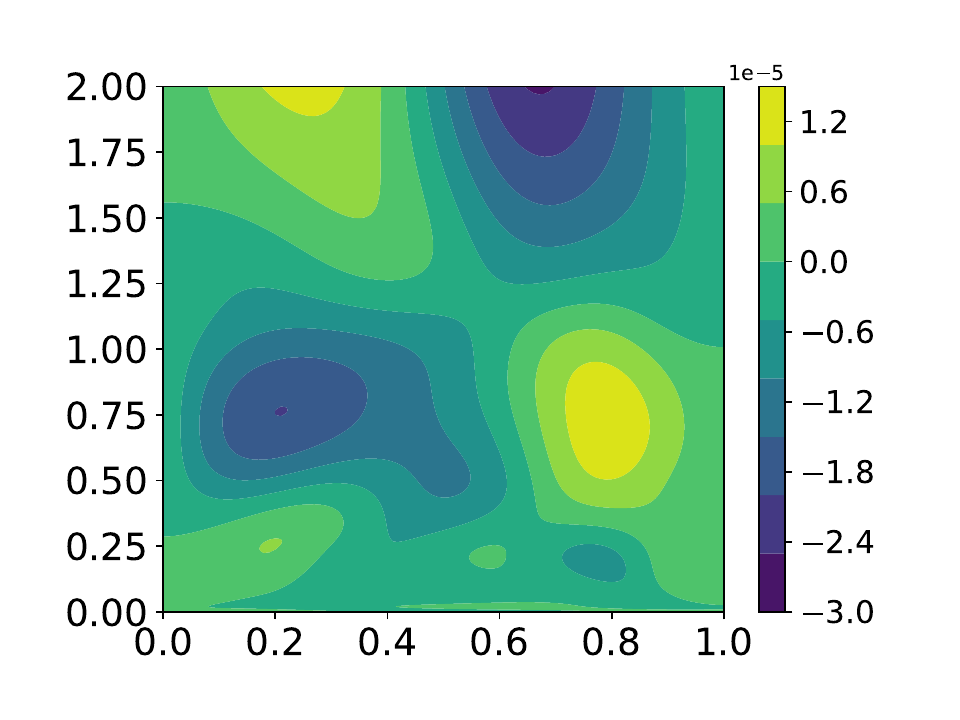}\hspace{-5mm}}
		\caption{Basis fields $g_i(\bm{x},t)$ of the DVS approximation.}
		\label{fig_burgers.4}
	\end{figure}

The basis fields $g_i(\bm{x},t)$ are shown in \Cref{fig_burgers.4}. Similarly, the core information of the solution $u(\bm{x},t;\bm{\xi})$ is captured by the first field $g_1(\bm{x},t)$, and the last few fields capture the fine-scale information of the original system.

\subsubsection{Allen-Cahn equation}
\label{num-ac}
At last, we consider the Allen-Cahn equation, which is defined on domain $D=[0,1]\times[0,1]$ as follows
	$$
		\left\{
		\begin{aligned}
			\frac{\partial u}{\partial t}&={\xi ^2} \Delta u -f(u),~~ \bm{x}\in D,~t\in[0,T],\\
			u(\bm{x},0;\bm{\xi})&=\sqrt{5}({x}_1^2-{x}_1)({x}_2^2-{x}_2),~~\bm{x}\in D,\\
			u(\bm{x},t;\bm{\xi})&=0,~~ \bm{x} \in \partial D,~t\in[0,T],
		\end{aligned}
		\right.
	$$
	with $T=1$, $f(u)=F^{'}(u)$, where $F(u)=\frac{1}{4}(u^2-1)^2$ is the Ginzburg-Landau double-well potential.
     We consider the setting that the parameter ${\xi}$ in the interval $[0.1,0.2]$, and the reference solution is calculated by the finite element method in space with mesh size $h_{{x}_1}=h_{{x}_2}=0.05$ and the backward Euler scheme in time with step size $\tau=10^{-4}$. 
     We choose $|\Xi|=8$ samples for training in \Cref{algorithm-dvs}.

\begin{figure}[htbp]
    \centering
    \includegraphics[width=0.5\textwidth, height=1.75in]{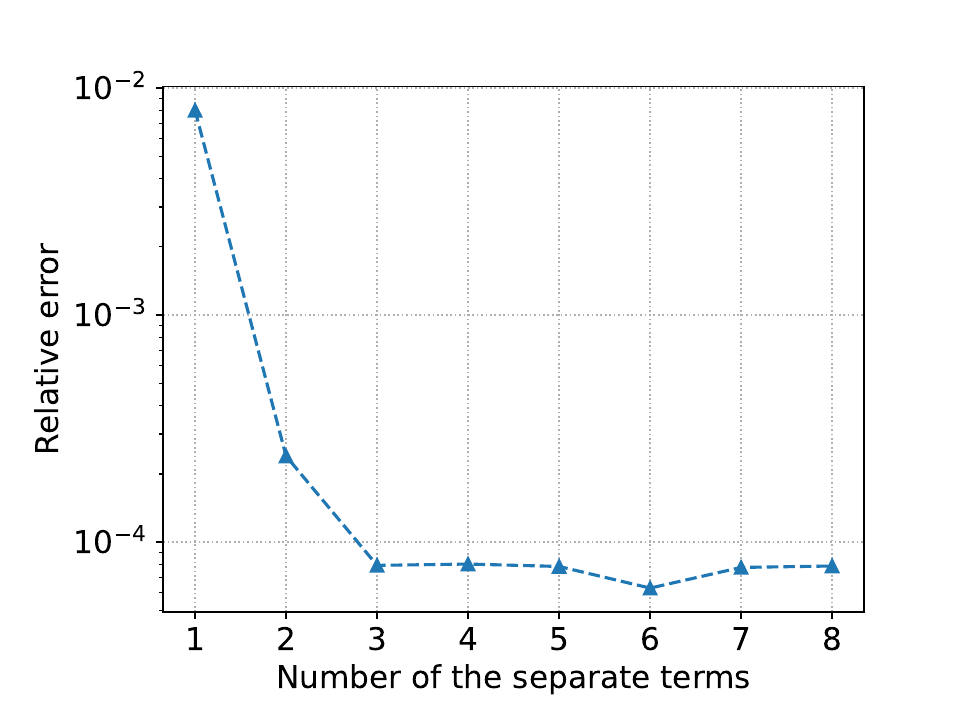}
    \caption{The average relative error corresponding to the different numbers of the separate terms $N$.}
    \label{fig_AC.1}
\end{figure}

First, we randomly choose $M=10^3$ samples and plot the average relative error calculated by equation (\ref{eq-meanerror}) versus the number of separate terms $N$ in \Cref{fig_AC.1}. The figure shows that the average relative error initially decreases with an increasing number of separate terms and then stabilizes.

\begin{figure}[htbp]
\centering
\subfigure[The probability density of the errors ]{
\includegraphics[width=2.4in, height=1.8in]{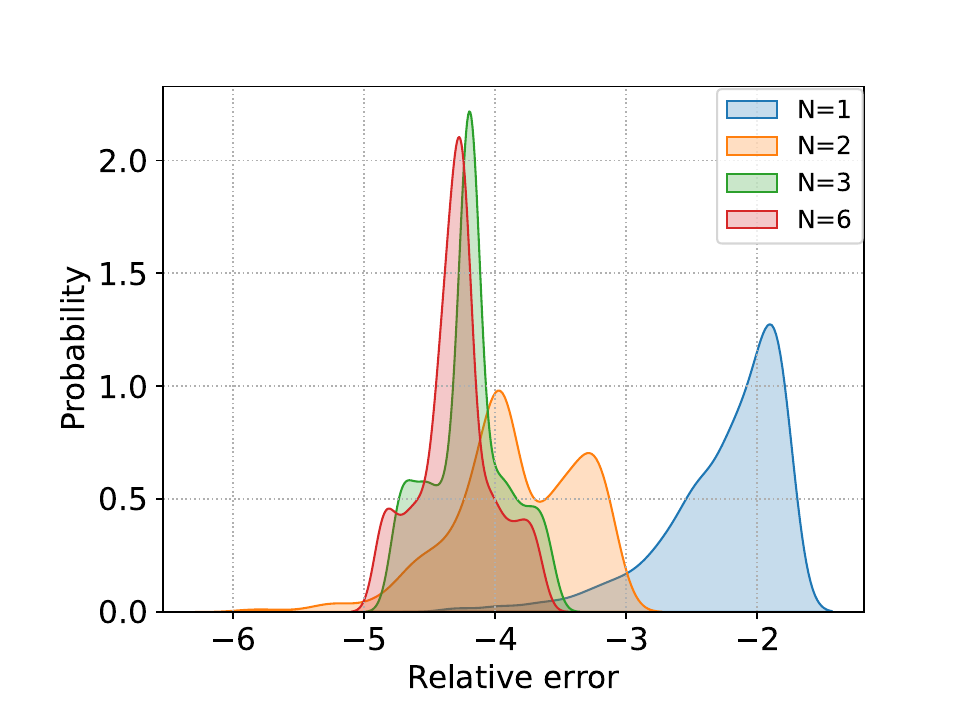}}
\subfigure[Error for the first 100 random samples]{
\includegraphics[width=2.4in, height=1.8in]{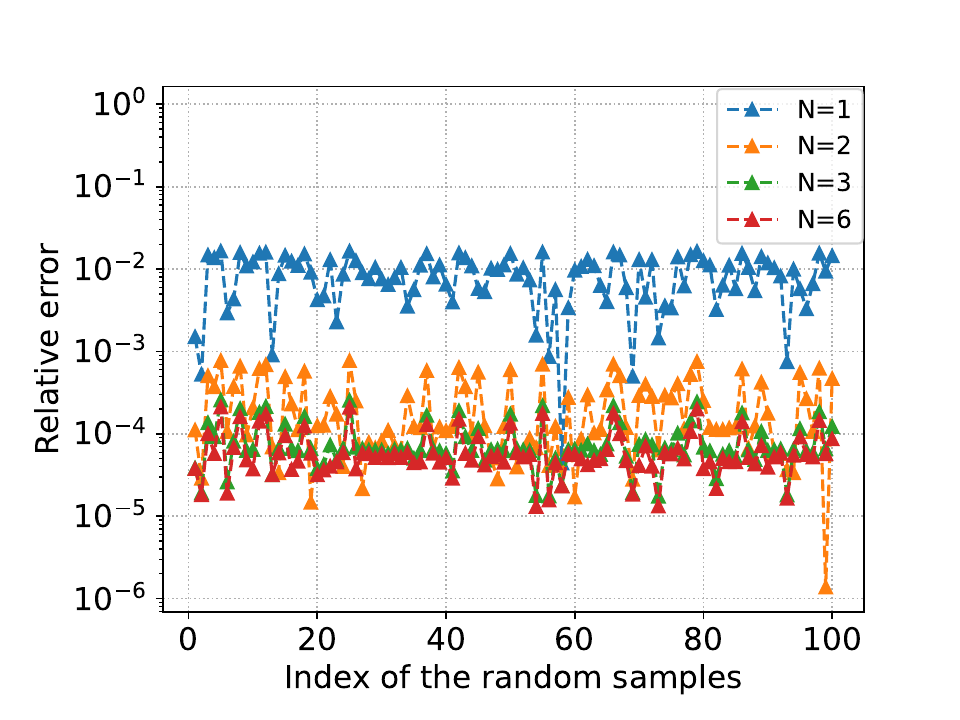}}
\caption{The relative error with the number of the separate terms $N=1$, $N=2$, $N=3$, and $N=6$.}
\label{fig_AC.2}
\end{figure}

Then we plot the distribution of the probability density of the logarithmic relative error (with base 10) for each random sample and the relative error for the first 100 parameter samples with the number of separate terms $N=1$, $N=2$, $N=3$, and $N=6$ in \Cref{fig_AC.2} to visualize the individual relative error. From the figure, we can see that: (1) the error is getting smaller and more compact as the number of separate terms $N$ increases; (2) the DVS method can give a good approximation for each sample.

\begin{figure}[htbp]
    \centering
    \includegraphics[width=.5\textwidth,height=1.8in]{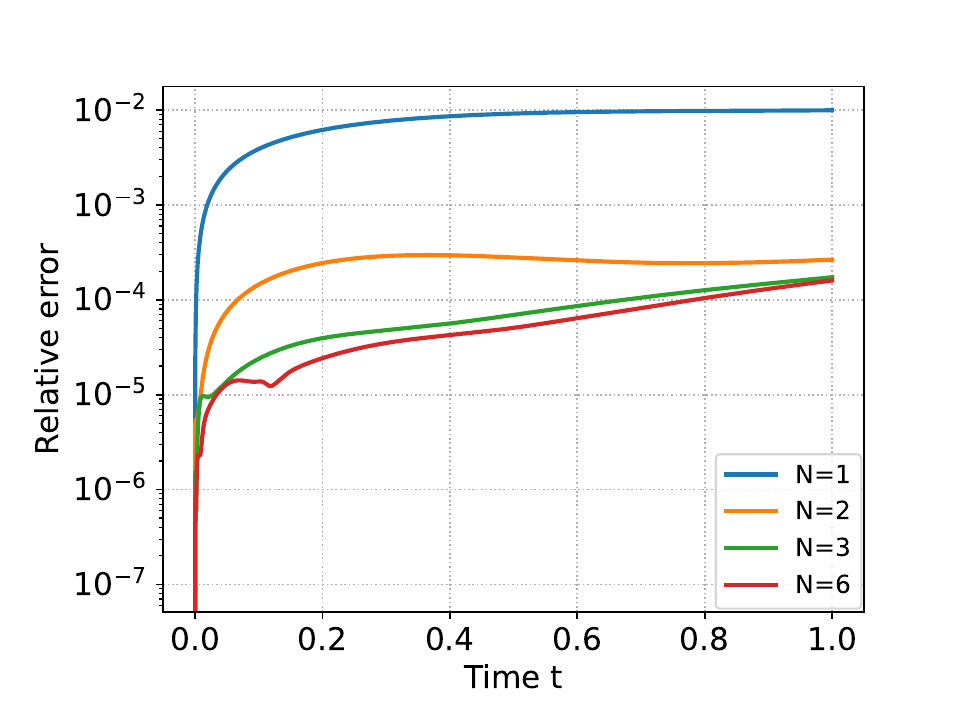}
    \caption{The average relative error versus different time levels with the number of the separate terms being $N=1$, $N=2$, $N=3$, and $N=6$.}
    \label{fig_AC.3}
\end{figure}

In \Cref{fig_AC.3}, we plot the average relative error versus different time levels with the number of separate terms being $N=1$, $N=2$, $N=3$, and $N=6$. It is shown that for each fixed time $t$, the error also decreases as the number of separate terms $N$ increases.

In \Cref{table_AC}, we list the average relative error and the average online CPU time based on $M=10^3$ random samples.
From the table, we can conclude that the approximation obtained by the DVS method can achieve a good trade-off in both approximation accuracy and computational efficiency for this numerical experiment.
\begin{table}[hbtp]
\centering
\caption{Comparison of the average relative errors and online CPU time for DVS and FEM-BE with different numbers of the separate terms.}
\vspace*{2pt}
\scriptsize
\begin{tabular}{cc|c|c}
\Xhline{1pt}			
\multicolumn{2}{c|}{\multirow{2}{*}\centering Algorithm}
&
\multicolumn{1}{c|}{\centering error $\epsilon$ }
&
\multicolumn{1}{c}{\centering online time}\\
\hline
\multirow{4}{*}{DVS}  & $N=1$  &   $7.99 \times10^{-3} $   &  $4.96 \times10^{-2}s$   \\
\cline{3-4}
                     & $N=2$    & $2.40 \times10^{-4} $   &  $9.39 \times10^{-2}s$  \\
\cline{3-4}
                     & $N=3$    & $7.91 \times10^{-5} $   &  $1.48 \times10^{-1}s$ \\
\cline{3-4}
                     & $N=4$    &   $8.02 \times10^{-5} $   &  $1.79 \times10^{-1}s$  \\
\hline
\multicolumn{2}{c|}{\centering FEM-BE} & $\setminus $ &  $26.54s$\\
\Xhline{1pt}
\end{tabular}
\label{table_AC}
\end{table}

\section{Conclusions}
\label{sec-Conclusions}
In this work, we proposed the dynamical Variable-separation (DVS) method for parameter-dependent dynamical systems to get a low-rank separate approximate representation of the solutions.
The proposed DVS method was devoted to constructing a separate approximation (\ref{eq-approx}) of the solution with a time-dependent spatial basis and a time-dependent parametric basis.
To obtain an efficient and reliable approximation, the DVS method contains two stages, the offline stage and the online stage.
In the offline stage, the proposed method adds the reduced basis functions through a greedy algorithm.
The method first constructs two uncoupled evolution equations directly derived from the original dynamical system and the previously separate representation terms, then modestly gives the initial conditions for the evolution equations at each enrichment step.
By solving the evolution equations, i.e., a parameter-independent PDE and a parameter-dependent ODE, the method can obtain the spatial basis and parameter-dependent basis,  respectively, which are all time-dependent.
Based on the separate approximation, the computational complexity of the online stage is independent of the spatial discretization of the original system, so the online stage is very efficient.
We applied the proposed DVS method to several parameter-dependent dynamical systems, and the numerical results suggest that the DVS method produces an efficient and robust reduced model.

Note that both the storage and the online computation of the DVS method rely on the size of time-space discretization. 
Reducing the cost of storage and computation requires further investigation. 
It is also of interest to see if the DVS method can tackle heterogeneous dynamical systems, where heterogeneity arises from physical properties, initial conditions, or other factors. 
Thus, it is also worth studying how to combine the DVS method with the partitioning of the spatial domain and/or the parameter space to handle more complex time-dependent nonlinear systems. 
Moreover, a rigorous convergence analysis of the proposed DVS method under reasonable assumptions is also of theoretical value. We leave these problems as our future research topics.

\appendix
\section{The definitions of \texorpdfstring{$\bigtriangleup_k(\bm{\xi})$}{}}
\label{sec-error}
An efficient error estimator $\bigtriangleup_k(\bm{\xi})$ is crucial for both the efficiency and the reliability of the proposed DVS method. In this part, we will describe two ways to define the error estimator $\bigtriangleup_k(\bm{\xi})$. 
The first one is very straightforward, i.e., one can take 
$\bigtriangleup_k(\bm{\xi}):=||e(\bm{x},t;\bm{\xi})||_{L^2([0,T];\mathcal{V})}$ with the error term given by 
$e(\bm{x},t;\bm{\xi}):=u(\bm{x},t;\bm{\xi})-u_{k-1}(\bm{x},t;\bm{\xi})$.
To further reduce the computational cost, we can follow \cite{Friess2017DynamicalMR,Hale1969pure, Wirtz2014a} to obtain posterior error bounds for the DVS method.
To this end, we first define the local logarithmic Lipschitz constant for the operator $\mathcal{F}$. 
\begin{definition}[Local logarithmic Lipschitz constant]
\label{def-local}
For a Lipschitz continuous function $\mathcal{F}:\mathcal{V}\to\mathcal{V^*}$, the local logarithmic Lipschitz constant of $\mathcal{F}$ at $u\in \mathcal{V}$ is defined by
$$
L_{\mathcal{V}}[\mathcal{F}](u): =
\mathop{\mathrm{sup}}\limits_{v\in \mathcal{V},v\neq u}\frac{\langle v-u,\mathcal{F}(v;\bm{\xi})-\mathcal{F}(u; \bm{\xi})\rangle}{||v-u||^2_\mathcal{V}}.
$$    
\end{definition}
 
Additionally, we should use the following comparison lemma, which is critical to derive the error estimator.
\begin{lemma}[Comparison lemma \cite{Wirtz2014a,Friess2017DynamicalMR}]
\label{lem-comparison}
Let $T>0$, and $u,\alpha,\beta:[0,T]\to R$ be integrable functions.
Suppose that $u$ is differentiable, and assume
$$
       \frac{\partial u}{\partial t}(t)\leq \beta(t) u(t)+\alpha(t),~~~\forall ~t\in[0,T].
$$
    Then one has 
$$
u(t)\leq\int_{0}^t\alpha(s)e^{\int_{s}^{t}\beta(\tau)\dd\tau} \dd s+e^{\int_{0}^{t}\beta(\tau)\dd\tau}u(0),~~\forall ~t\in [0,T].
$$
\end{lemma}

From the definition of the parameter-dependent dynamical system \eqref{eq-dynamical-system} we know that the approximate error $e(\bm{x},t;\bm{\xi})$ satisfies
\begin{equation}
\label{eq-error-estimate}
\frac{\partial e}{\partial t}(\bm{x},t;\bm{\xi}) =\mathcal{F}(u(\bm{x},t;\bm{\xi});\bm{\xi})-\mathcal{F}(u_{k-1}(\bm{x},t;\bm{\xi});\bm{\xi})+r_k(\bm{x},t;\bm{\xi}),
\end{equation}
where $r_k(\bm{x},t;\bm{\xi})=\mathcal{F}(u_{k-1}(\bm{x},t;\bm{\xi});\bm{\xi})-\frac{\partial u_{k-1}}{\partial t}(\bm{x},t;\bm{\xi})$.
Based on \Cref{lem-comparison}, by combining equation (\ref{eq-error-estimate}) with the definition of the local logarithmic Lipschitz constant for $\mathcal{F}$, we can derive a bound of the error ${e}(\bm{x},t;\bm{\xi})$ as follows.
\begin{proposition}
\label{prop-e}
The error norm $||{e}(\bm{x},t;\bm{\xi}) ||_{\mathcal{V}}$ satisfies 
$$
||{e}(\bm{x},t;\bm{\xi})||_{\mathcal{V}}\leq\delta_k(t;\bm{\xi}),
$$
\end{proposition}
where 
\begin{equation}
\label{errordelta}
\delta_k(t;\bm{\xi})=\int_{0}^t\alpha(s;\bm{\xi})e^{\int_{s}^{t}\beta(\tau;\bm{\xi})\dd\tau}\dd s+e^{\int_{0}^{t}\beta(\tau;\bm{\xi})\dd\tau}\|e(\bm{x},0;\bm{\xi})\|_{\mathcal{V}},
\end{equation}
and $\alpha(t;\bm{\xi})=||{r}_k(\bm{x},t;\bm{\xi})||_{\mathcal{V}}$, $\beta(t;\bm{\xi})=L_{\mathcal{V}}[\mathcal{F}]({u}_{k-1}(\bm{x},t;\bm{\xi}))$.

\begin{proof}
By the definition of the approximate error ${e}(\bm{x},t;\bm{\xi})$, we have
$$
\begin{aligned}
  \frac{1}{2}\frac{\dd}{\dd t}||{e}(\bm{x},t;\bm{\xi})||_{\mathcal{V}}^2&=\langle {e}(\bm{x},t;\bm{\xi}),\frac{\partial e}{\partial t}(\bm{x},t;\bm{\xi}) \rangle \\
  &=\langle {e}(\bm{x},t;\bm{\xi}), \mathcal{F}({u}(\bm{x},t;\bm{\xi});\bm{\xi})-\mathcal{F}({u}_{k-1}(\bm{x},t;\bm{\xi});\bm{\xi})+{r}_k(\bm{x},t;\bm{\xi})\rangle\\
  &\leq L_\mathcal{V}[\mathcal{F}]({u}_{k-1}(\bm{x},t;\bm{\xi}))||{e}(\bm{x},t;\bm{\xi})||_{\mathcal{V}}^2+||{r}_k(\bm{x},t;\bm{\xi})||_{\mathcal{V}}||{e}(\bm{x},t;\bm{\xi})||_{\mathcal{V}},
\end{aligned}
$$
where the last inequality is due to the definition of the local logarithmic Lipschitz constant of $\mathcal{F}$ in \Cref{def-local}. Then we can obtain that
$$
    \frac{\dd}{\dd t}||{e}(\bm{x},t;\bm{\xi})||_{\mathcal{V}}\leq L_\mathcal{V}[\mathcal{F}]({u}_{k-1}(\bm{x},t;\bm{\xi}))||{e}(\bm{x},t;\bm{\xi})||_{\mathcal{V}}+||{r}_k(\bm{x},t;\bm{\xi})||_{\mathcal{V}}.
$$
According to \Cref{lem-comparison}, we finish the proof.  
\end{proof}

Thanks to \Cref{prop-e}, one can take
$\bigtriangleup_k(\bm{\xi}):=\sqrt{\int_{0}^{T}(\delta_k(t;\bm{\xi}))^2\dd t}$ as the error estimator with $\delta_k(t;\bm{\xi})$ given by \eqref{errordelta}.

\begin{remark}
To efficiently calculate the error estimator, one can further use an offline-online procedure. 
The details can be found in \cite{Friess2017DynamicalMR}.
\end{remark}

    \section{The Variable-separation method for heat equation}
    \label{sec-vs}
    Here we discuss the details of the Variable-separation (VS) method for the heat equation (presented in  \cref{num-heat}) to obtain the separate approximation
$$
u(\bm{x},t;\bm{\xi})\approx u_N(\bm{x},t;\bm{\xi}):=\sum_{i=1}^{N}\zeta_i(\bm{\xi})g_i(\bm{x},t),
$$
    where each $\zeta_i(\bm{\xi})$ is a parameter-dependent but time-independent basis function, each $g_i(\bm{x},t)$ is a time-dependent spatial basis function.
    In this setting, the variational formulation \eqref{eq-varia} becomes
    \begin{equation}
    \label{eq-var-heat}
     \Big\langle \frac{\partial u(\bm{x},t;\bm{\xi})}{\partial t},v\Big\rangle
        +a(u(\bm{x},t;\bm{\xi}),v;\bm{\xi})=\langle f(\bm{x};\bm{\xi}), v\rangle, ~\forall ~v\in \mathcal{V},~ t\in [0,T],
    \end{equation}
    where $a(u,v;\bm{\xi})=-\kappa(\bm{\xi})\langle \Delta u, v\rangle$, which is affine with respect to $\bm{\xi}$.
    
    Define the error of the approximation as ${e}(\bm{x},t;\bm{\xi}):={u}(\bm{x},t;\bm{\xi})-{u}_{k-1}(\bm{x},t;\bm{\xi})$ with $u_{k-1} \equiv 0$ for $k=1$, and the residual
$$
{r}_k(\bm{x},t;\bm{\xi}):=
\begin{cases}
f(\bm{x};\bm{\xi}),&k=1,\\
f(\bm{x};\bm{\xi}) + \kappa(\bm{\xi})\Delta u_{k-1}(\bm{x},t;\bm{\xi})-\frac{\partial u_{k-1}}{\partial t}(\bm{x},t;\bm{\xi}),&k\geq2.
\end{cases}
$$
Then, equation (\ref{eq-var-heat}) turns to
\begin{equation}
\label{eq-residual-2}
\Big \langle \frac{\partial e}{\partial t}(\bm{x},t;\bm{\xi}),v\Big\rangle+a(e(\bm{x},t;\bm{\xi}),v;\bm{\xi})=\langle{r}_k(\bm{x},t;\bm{\xi}),v\rangle.
\end{equation}
  
At step $k$, we choose $\bm{\xi}_k\in {\text{argmax}}_{\bm{\xi}\in\Xi} \bigtriangleup_k(\bm{\xi})$ following the strategy in \cref{sec-dvs-framework}, and take $g_k(\bm{x},t)$  as the solution of (\ref{eq-residual-2}) with $\bm{\xi}=\bm{\xi}_k$, subject to the initial condition
$$
g_k(\bm{x},0):= \begin{cases}
\begin{aligned}
&\mu(\bm{x};\bm{\xi}_1), &k=1,\\
&\mu(\bm{x};\bm{\xi}_k)-\sum_{i=1}^{k-1}\zeta_i(\bm{\xi}_k)g_{i}(\bm{x},0), &k\geq 2.
\end{aligned}
\end{cases}$$ 
Let $\zeta_k(\bm{\xi})$ be an unknown basis function and define $\tilde{e}(\bm{x},t;\bm{\xi}):=g_k(\bm{x},t)\zeta_k(\bm{\xi})$. Taking $e=\tilde{e}$ and $v=g_k(\bm{x},t)$ in equation (\ref{eq-residual-2}), we have
    \begin{equation}
		\label{eq-e-xi-2}
			\quad \Big(\Big\langle \frac{\partial g_k(\bm{x},t)}{\partial t},g_k(\bm{x},t)\Big\rangle +a(g_k(\bm{x},t),g_k(\bm{x},t);\bm{\xi})\Big)\zeta_k(\bm{\xi})=\langle{r}_k(\bm{x},t;\bm{\xi}),g_k(\bm{x},t)\rangle.
	\end{equation}
In the practical simulation, by randomly choosing $n_t$ time samples (where $n_t\ll N_t$, $N_t$ is the number of the time discretization used to calculate (\ref{eq-residual-2})), we can derive the explicit expression of $\zeta_k(\bm{\xi})$ corresponding to each time sample from (\ref{eq-e-xi-2}).
Here, we choose the basis $\zeta_k(\bm{\xi})$ which makes $\bigtriangleup_{k+1}(\bm{\xi})$ the smallest.

\bibliographystyle{plain}

\begin{thebibliography}{}

    \bibitem{C2013nonintrusive}
    {\sc C. Audouze, F. De Vuyst, and P.B. Nair}, {\em Nonintrusive reduced-order modeling of parametrized time-dependent partial differential equations}, Numer. Methods Partial Differential Equations, 29 (2013), pp. 1587-1628.
    
    
    \bibitem{U2014model}
    {\sc U. Baur, P. Benner, and L. Feng}, {\em Model order reduction for linear and nonlinear systems: a system-theoretic perspective}, Arch. Comput. Methods Eng., 21 (2014), pp. 331–358. 


    \bibitem{Beck2000the}
    {\sc M.H. Beck, A. J{\"a}ckle, G.A. Worth, and H.-D. Meyer}, {\em The multiconfiguration time-dependent Hartree (MCTDH) method: a highly efficient algorithm for propagating wavepackets}, Phys. Rep., 324 (2000), pp. 1-105.
    
    \bibitem{P2015a}
    {\sc P. Benner, S. Gugercin, and K. Willcox}, {\em A survey of projection-based model reduction methods for parametric dynamical systems}, SIAM Rev., 57 (2015), pp. 483-531. 

    \bibitem{Friess2017DynamicalMR}
    {\sc M. Billaud-Friess and A. Nouy}, {\em Dynamical model reduction method for solving parameter-dependent dynamical systems}, SIAM J. Sci. Comput., 39 (2017), pp. A1766-A1792.
    
  
    \bibitem{T2008model}
    {\sc T. Bui-Thanh, K. Willcox, and O. Ghattas}, {\em Model reduction for large-scale systems with high-dimensional parametric input space}, SIAM J. Sci. Comput., 30 (2008), pp. 3270–3288.

    
    \bibitem{Chen2024stochastic}
    {\sc L. Chen, Y. Chen, Q. Li, and Z. Zhang}, {\em Stochastic domain decomposition based on variable-separation method}, Comput. Methods Appl. Mech. Engrg., 418 (2024), 116538.
    
    
    \bibitem{Cheng2013a1}
    {\sc M. Cheng, T. Y. Hou, and Z. Zhang}, {\em A dynamically bi-orthogonal method for time-dependent stochastic partial differential equations I: Derivation and algorithms}, J. Comput. Phys., 242 (2013), pp. 843-868.
    
    \bibitem{Cheng2013a2}
    {\sc M. Cheng, T. Y. Hou, and Z. Zhang}, {\em A dynamically bi-orthogonal method for time-dependent stochastic partial differential equations II: Adaptivity and generalizations}, J. Comput. Phys., 242 (2013), pp. 753-776.

    \bibitem{L2006Semi‐implicit}
    {\sc L. Diening, A. Prohl, and M. Růžička}, {\em Semi‐implicit Euler scheme for generalized Newtonian fluids}, SIAM J. Numer. Anal., 44 (2006), pp. 1172–1190.

    \bibitem{Feppon2018a}
    {\sc F. Feppon and P.F.J. Lermusiaux}, {\em A geometric approach to dynamical model order reduction}, SIAM J. Matrix Anal. Appl., 39 (2018), pp. 510-538.

      
    \bibitem{R2019reduced}
    {\sc R. Ghaffari and F. Ghoreishi}, {\em Reduced collocation method for time-dependent parametrized partial differential equations}, Bull. Iran. Math. Soc., 45 (2019), pp. 1487–1504.

     \bibitem{C2019decay}
    {\sc C. Greif and K. Urban}, {\em Decay of the Kolmogorov $N$-width for wave problems}, Appl. Math. Lett., 96 (2019), pp. 216-222.
    
     \bibitem{M2019data}
    {\sc M. Guo and J.S. Hesthaven}, {\em Data-driven reduced order modeling for time-dependent problems}, Comput. Methods Appl. Mech. Engrg., 345 (2019), pp. 75-99.


     \bibitem{Hale1969pure}
    {\sc J.K. Hale}, {\em Ordinary Diﬀerential Equations}, Pure Appl. Math. XXI, Wiley-Interscience, New York, 1969.

   
    \bibitem{Hesthaven2022reduced}
    {\sc J.S. Hesthaven, C. Pagliantini, and G. Rozza}, {\em Reduced basis methods for time-dependent problems}, Acta Numer., 31 (2022), pp. 265-345.

    \bibitem{V2020analy}
    {\sc V.H. Hoang}, {\em Analyticity, Regularity, and generalized polynomial chaos approximation of stochastic, parametric parabolic two-scale partial differential equations}, Acta Math. Vietnam., 45 (2020), pp. 217–247.
    
    \bibitem{Hoang2012REGULARITYAG}
    {\sc V.H. Hoang and C. Schwab}, {\em  Regularity and generalized polynomial chaos approximation of parametric and random second-order hyperbolic partial differential equations}, Anal. Appl. (Singap.), 10 (2012), pp. 295-326.


     \bibitem{Jiang2018model}
    {\sc L. Jiang and Q. Li}, {\em Model reduction method using variable-separation for stochastic saddle point problems}, J. Comput. Phys., 354 (2018), pp. 43-66.
    
    \bibitem{Koch2007regularity}
    {\sc O. Koch and C. Lubich}, {\em Regularity of the multi-configuration time-dependent Hartree approximation in quantum molecular dynamics}, ESAIM, Math. model. numer. anal., 41 (2007), pp. 315–331.

     \bibitem{k2020model}
    {\sc K. Lee and K.T. Carlberg}, {\em Model reduction of dynamical systems on nonlinear manifolds using deep convolutional autoencoders}, J. Comput. Phys., 404 (2020), 108973.


     \bibitem{Li2017a}
    {\sc Q. Li and L. Jiang}, {\em A novel variable-separation method based on sparse and low rank representation for stochastic partial differential equations}, SIAM J. Sci. Comput., 39 (2017), pp. A2879-A2910.

     \bibitem{Li2024an}
    {\sc Q. Li, C. Liu, M. Li, and P. Zhang}, {\em An adaptive method based on local dynamic mode decomposition for parametric dynamical systems}, Commun. Comput. Phys., 35 (2024), pp. 38-69.
    
    \bibitem{Li2020a}
    {\sc Q. Li and P. Zhang}, {\em A variable-separation method for nonlinear partial differential equations with random inputs}, SIAM J. Sci. Comput., 42 (2020), pp. A723-A750.

     
    
    \bibitem{Musharbash2018dual}
    {\sc E. Musharbash and F. Nobile}, {\em Dual dynamically orthogonal approximation of incompressible Navier Stokes equations with random boundary conditions}, J. Comput. Phys., 354 (2018), pp. 135-162.

     \bibitem{Musharbash2015error}
    {\sc E. Musharbash, F. Nobile, and T. Zhou}, {\em Error analysis of the dynamically orthogonal approximation of time dependent random PDEs}, SIAM J. Sci. Comput., 37 (2015), pp. A776-A810.
    
     \bibitem{K2020long}
    {\sc K. Nath, A. Dutta, and B. Hazra}, {\em Long duration response evaluation of linear structural system with random system properties using time dependent polynomial chaos}, J. Comput. Phys., 418 (2020), 109596.

     \bibitem{Nouy2009recent}
    {\sc A. Nouy}, {\em Recent developments in spectral stochastic methods for the numerical solution of stochastic partial differential equations}, Arch. Comput. Methods Eng., 16 (2009), pp. 251–285.

    \bibitem{Nouy2009gener}
    {\sc A. Nouy and O.P.Le Ma{\^i}tre}, {\em Generalized spectral decomposition for stochastic nonlinear problems}, J. Comput. Phys., 228 (2009), pp. 202-235.
    
    \bibitem{M2016reduced}
    {\sc M. Ohlberger and S. Rave}, {\em Reduced basis methods: success, limitations and future challenges}, Proceedings of the Conference Algoritmy, (2016), pp. 1-12.
  
    \bibitem{Othmar2007dynamical}
    {\sc K. Othmar and L. Christian}, {\em Dynamical low-rank approximation}, SIAM J. Matrix Anal. Appl., 29 (2007), pp. 434-454.


   
    \bibitem{Patil2020real-time}
    {\sc P. Patil and H. Babaee}, {\em Real-time reduced-order modeling of stochastic partial differential equations via time-dependent subspaces}, J. Comput. Phys., 415 (2020), 109511.

    \bibitem{B2022breaking}
    {\sc B. Peherstorfer}, {\em Breaking the Kolmogorov barrier with nonlinear model reduction}, Notices Amer. Math. Soc., 69 (2022), pp. 725-733.


    \bibitem{A2012N-Widths}
    {\sc A. Pinkus}, {\em N-widths in approximation theory}, Springer Berlin, Heidelberg, 2012.


    \bibitem{D2023manifold}
    {\sc D. Rim, B. Peherstorfer, and K.T. Mandli}, {\em Manifold approximations via transported subspaces: model reduction for transport-dominated problems}, SIAM J. Sci. Comput., 45 (2023), pp. A170-A199.
    
    \bibitem{E2008algebraic}
    {\sc E. Rosseel, T. Boonen, and S. Vandewalle}, {\em  Algebraic multigrid for stationary and time-dependent partial differential equations with stochastic coefficients}, Numer. Linear Algebra Appl., 15 (2008), pp. 141-163.

     \bibitem{Sapsis2009dynamical}
    {\sc T.P. Sapsis and P.F.J. Lermusiaux}, {\em Dynamically orthogonal field equations for continuous stochastic dynamical systems}, Phys. D, 238 (2009), pp. 2347-2360.
    
    \bibitem{Sapsis2012dynamical}
    {\sc T.P. Sapsis and P.F.J. Lermusiaux}, {\em Dynamical criteria for the evolution of the stochastic dimensionality in flows with uncertainty}, Phys. D, 241 (2012), pp. 60-76.

    \bibitem{Sleeman2022goal}
    {\sc M.K. Sleeman and M. Yano}, {\em Goal-oriented model reduction for parametrized time-dependent nonlinear partial differential equations}, Comput. Methods Appl. Mech. Engrg., 388 (2022), 114206.
    
    \bibitem{Song2024a}
    {\sc H. Song, Y. Ba, D. Chen, and Q. Li}, {\em  A model reduction method for parametric dynamical systems defined on complex geometries}, J. Comput. Phys., 506 (2024), 112923.

     \bibitem{X2021non}
    {\sc X. Sun and J. Choi}, {\em Non-intrusive reduced-order modeling for uncertainty quantification of space–time-dependent parameterized problems}, Comput. Math. Appl., 87 (2021), pp. 50-64.

    
    \bibitem{Ueckermann2013numerical}
    {\sc M.P. Ueckermann, P.F.J. Lermusiaux, and T.P. Sapsis}, {\em Numerical schemes for dynamically orthogonal equations of stochastic fluid and ocean flows}, J. Comput. Phys., 233 (2013), pp. 272-294.

    
    \bibitem{B2019Kolmo}
    {\sc B. Unger and S. Gugercin}, {\em Kolmogorov n-widths for linear dynamical systems}, Adv. Comput. Math., 45 (2019), pp. 2273–2286.
    
    \bibitem{Wirtz2014a}
    {\sc  D. Wirtz, D.C. Sorensen, and B. Haasdonk}, {\em A posteriori error Estimation for DEIM reduced nonlinear dynamical systems}, SIAM J. Sci. Comput., 36 (2014), pp. A311–A338.
    
     \bibitem{Xiu2002the}
    {\sc D. Xiu and G.E. Karniadakis}, {\em  The Wiener--Askey polynomial chaos for stochastic differential equations}, SIAM J. Sci. Comput., 24 (2002), pp. 619-644.

    \bibitem{Xiu2003modeling}
    {\sc D. Xiu and G.E. Karniadakis}, {\em  Modeling uncertainty in flow simulations via generalized polynomial chaos}, J. Comput. Phys., 187 (2003), pp. 137-167.

    
   
\end{thebibliography}

\end{document}